\documentclass{daj}
\usepackage{amssymb}
\usepackage{latexsym}
\usepackage{mathrsfs}

\newtheorem{theorem}{Theorem}[section]
\newtheorem{definition}[theorem]{Definition}
\newtheorem{proposition}[theorem]{Proposition}
\newtheorem{lemma}[theorem]{Lemma}
\newtheorem{corollary}[theorem]{Corollary}
\newtheorem{remark}[theorem]{Remark}

\newtheorem{question}[theorem]{Question}
\newtheorem{property}[theorem]{Property}
\newcommand {\n}{{\mathbb N}}
\newcommand {\V}{{\mathbb V}}

\newcommand {\real}{{\mathbb R}}

\newcommand {\Q}{{\mathcal{Q}}}

\newcommand {\p}{{\mathcal{P}}}

\newcommand {\R}{{\mathcal{R}}}
\newcommand {\SD}{S\!D}
\newcommand {\res}{\!\upharpoonright\!}
\newcommand {\f}{{\mathcal{F}}}

\dajAUTHORdetails{%
  title = {A Simple Combinatorial Proof of Szemer\'{e}di's Theorem
via Three Levels of Infinities}, 
  author = {Renling Jin},
  plaintextauthor = {Renling Jin},
  plaintexttitle = {A Simple Combinatorial Proof of Szemer\'{e}di's Theorem
via Three Levels of Infinities}, 
  runningtitle = {A Simple Proof},
  runningauthor = {Renling Jin},
  copyrightauthor = {Renling Jin},
  keywords = {arithmetic progression, Szemer\'{e}di's theorem,
nonstandard analysis, iterated nonstandard extensions},
}   

\dajEDITORdetails{%
   year={2023},
   number={15},
   received={11 April 2022},   
   revised={5 July 2023},    
   published={25 September 2023},  
   doi={10.19086/da.87772},       
}   

\begin{document}

\begin{frontmatter}[classification=text]

\title{A Simple Combinatorial Proof of Szemer\'{e}di's Theorem
via Three Levels of Infinities} 

\author[pgom]{Renling Jin\thanks{This work was partially supported by
a collaboration grant (ID: 513023) from Simons Foundation.}}

\begin{abstract}
We present a nonstandard simple elementary proof of
Szemer\'{e}di's theorem by a straightforward induction
with the help of three levels of infinities and four
different bounded elementary embeddings in a nonstandard universe.
\end{abstract}
\end{frontmatter}


\section{Introduction}
 \quad This article is strongly influenced by Terence Tao's notes \cite{tao}.

\begin{theorem}[van der Waerden, 1927]
Given any $k,n\in\n$, there exists
$\Gamma(k,n)\in\n$ called van der Waerden number,
such that if $\{U_1,U_2,\ldots,U_n\}$ is a partition
of $\{1,2,\ldots,\Gamma(k,n)\}$, then there is a $u\leq n$ such that
$U_u$ contains a $k$--term arithmetic progression.
\end{theorem}

\begin{theorem}[E.\ Szemer\'{e}di, 1975 \cite{szemeredi}]
If $D\subseteq\n$ has a positive upper density, then
$D$ contains a $k$--term arithmetic progression for every $k\in\n$.
\end{theorem}

Szemer\'{e}di's theorem confirms a conjecture of P.\ Erd\H{o}s
and P.\ Tur\'{a}n made in 1936, which implies van der Waerden's theorem.

Nonstandard versions of Furstenberg's ergodic proof and Gowers's
harmonic proof of Szemer\'{e}di's theorem have been tried by T.\ Tao
(see Tao's blog post \cite{tao2}).
In the workshop {\it Nonstandard methods in
combinatorial number theory} sponsored by
American Institute of Mathematics in San Jose, CA,
 August 2017, Tao gave a series of lectures to explain
Szemer\'{e}di's original combinatorial proof and hope
to simplify it so that the proof can be
better understood. He believed that Szemer\'{e}di's combinatorial
method should have a greater impact on combinatorics.

During these lectures Tao challenged the audience
to produce a nonstandard proof of Szemer\'{e}di's theorem which is
noticeably simpler and more transparent than Szemer\'{e}di's
original proof. The current article is the result of Tao's challenge
and inspiration. However, in his later blog post \cite{tao3},
Tao commented that ``in fact there are now signs that perhaps nonstandard
analysis is not the optimal framework in which to place this argument.''
We disagree. The current article is our effort to show that
with the help of a nonstandard universe with three levels
of infinities, Szemer\'{e}di's original argument
can be made simpler and more transparent.

The main simplification in our proof of Szemer\'{e}di's theorem
compared to the standard proof in
\cite{szemeredi,tao} is that a Tower of Hanoi type induction
in \cite[Theorem 6.6]{tao} and in \cite[Lemma 5, Lemma 6, and Fact 12]{szemeredi}
is replaced by a straightforward induction (see Lemma \ref{key} below),
which makes Szemer\'{e}di's idea more transparent.
To achieve this, we work within a chain of nonstandard extensions
$\V_0\prec\V_1\prec\V_2\prec\V_3$ which supply three levels of
infinities, plus various bounded elementary embeddings from
$\V_j$ to $\V_{j'}$ for some $0\leq j<j'\leq 3$.

The paper is organized in the following sections. \S\ref{nsa}
is a brief introduction of logic foundation for constructing
the nonstandard extensions $\V_0\prec\V_1\prec\V_2\prec\V_3$ and
bounded elementary embeddings including
$i_0$, $i_*$, $i_1$, and $i_2$ to the reader who does not have
logic background. The reader who is only interested in
applications can familiarize with the notation,
Property \ref{property}, Proposition \ref{internal2},
and Proposition \ref{spill2}, safely skip the proof,
and return to it at a later time.
In \S\ref{nsatools} we translate the density along arithmetic progressions
in standard setting to strong upper Banach density in nonstandard setting
as well as translate some consequences of the
double counting argument in standard setting to nonstandard setting.
The reader who is only interested in
applications can familiarize with the notation,
Lemma \ref{equiv}, Lemma \ref{bigD}, and Lemma \ref{bigS},
safely skip the proof, and return to it at a later time.
In \S\ref{RegularitySection} we re-write,
in a nonstandard setting, the proof of a so--called
mixing lemma in \cite{tao} based on a weak regularity lemma.
This section does not offer new idea but is included only for self-containment.
\S\ref{allk} is the main part of the paper where
we present the proof of Szemer\'{e}di's theorem. In \S\ref{question}
we pose a question whether the presented proof of Szmer\'{e}di's theorem
can be carried out without the axiom of choice.

\section{Construction of Nonstandard Extensions}\label{nsa}

\quad The notation we use here should be consistent with some standard
 textbooks. Consult, for example, \cite{CK,goldblatt,LW}
for more details. If $f:A\to B$ is a function, then $f(a)$ denotes
the image of $a$ as an element in $B$
and $f[C]:=\{f(a)\mid a\in C\}$ for some $C\subseteq A$.

\bigskip

\noindent {\bf \S\ref{nsa}.1 Superstructure}\quad
Let $\omega$ be the set of all standard non-negative integers used
at the meta-level and $X$ be an infinite set of urelements,
i.e., elements without members. The {\it superstructure} on $X$
(cf.\ \cite[page 263]{CK}), denoted by
$\V(X)$, is composed of the base set $V(X)$ and the membership relation
$\in$ on $V(X)$ where
\[V(X):=\bigcup_{n\in\omega}V(X,n)\] and
$V(X,n)$ is defined recursively by letting
\[V(X,0):=X\,\mbox{ and }\,V(X,n+1):=V(X,n)\cup\mathscr{P}(V(X,n))\]
for every $n\in\omega$ where $\mathscr{P}$ is the powerset operator.
For notational convenience we often write $\V(X)$ also for its base set $V(X)$.
Hopefully, this will not cause confusion. One can define a rank function for
every $a\in\V(X)$ recursively. Set $\mbox{rank}(a)=0$ iff (the abbreviation of
``if and only if'') $a\in X$ and
\begin{equation}\label{rank}
\mbox{rank}(a)=1+\max\{\mbox{rank}(b)\mid b\in a\}
\end{equation} for any $a\in\V(X)\setminus X$.
Notice that $\mbox{rank}(a)=n$ iff $a\in \V(X,n)\setminus \V(X,n-1)$.

Let $\n_0$ be the set of all standard positive
integers and $\real_0$ be the set of all standard real numbers.
By the {\it standard universe} we mean the superstructure $\V_0:=\V(\real_0)$
on $X=\real_0$.
Notice that all standard mathematical objects mentioned in this paper
have ranks below, say, $100$. For example, an ordered pair $(a,b)$ of
standard real numbers can be viewed as the set $\{\{a\},\{a,b\}\}\in\V(\real_0,2)$
and a function $f:\n_0\to\real_0$ can be viewed as a set of ordered
pairs in $\V(\real_0,2)$. Hence $f\in\V(\real_0,3)$.

\bigskip

\noindent {\bf \S\ref{nsa}.2 Logic and model theory}\quad
Before introducing nonstandard universe we should mention
briefly, without rigor, some concepts in model theory.
For simplicity we consider only model theory
on finite relational languages. We call a set $\mathscr{L}$ of
finitely many symbols $P_1,P_2,\ldots,P_n$ with arity $m_i\in\omega$
for each $P_i$ a (relational) language.
An $\mathscr{L}$--model $\mathcal{M}$ is a structure composed of
a nonempty base set $M$ and an $m_i$-nary relation
$P^{\mathcal{M}}_i\subseteq M^{m_i}$, called the interpretation of
$P_i$ in $\mathcal{M}$, for
each symbol $P_i\in\mathscr{L}$. For notational convenience
we sometimes write $\mathcal{M}$ for
the base set of $\mathcal{M}$.

We can define (first-order)
$\mathscr{L}$--formulas recursively starting from atomic formulas. If
$x$ and $x'$ are variables, then
$x=x'$ is an atomic formula.
If $P_i\in\mathscr{L}$ and $\overline{x}:=(x_1,x_2,\ldots,x_{m_i})$ is
an $m_i$-tuple of variables, then $P_i(\overline{x})$ is an atomic formula.
The word ``first-order'' means that these variables are intended to
take only elements of some $\mathscr{L}$--model as their values. All formulas mentioned in this paper are first-order.
We will use the symbol $\overline{a}$ to represent
an $m$-tuple of elements with some suitable generic number $m\in\omega$.
When $\overline{x}$ is intended to be substituted by $\overline{a}$,
we assume implicitly that they have the same length.

If $\varphi$
and $\psi$ are $\mathscr{L}$--formulas and $x$ is a variable,
then $\neg\varphi$, $\varphi\wedge\psi$, $\varphi\vee\psi$,
$\varphi\to\psi$,  $\varphi\leftrightarrow\psi$, $\forall x\varphi$, and
$\exists x\varphi$ are also $\mathscr{L}$--formulas. The symbols
$\neg$, $\wedge$, $\vee$, $\to$, and $\leftrightarrow$ are called
logic connectives, and $\forall$ and $\exists$ are called universal and
existential quantifiers. An occurrence of
$x$ in $\varphi$ is called {\it bounded} if it occurs in a sub-formula of the form
$\forall x\psi$ or $\exists x\psi$ in $\varphi$.
An occurrence of $x$ in $\varphi$ is called free if it is not bounded.
A formula without free variable is called a sentence.

Let $\mathcal{M}$ be an $\mathscr{L}$--model. For any $\mathscr{L}$--formula
$\varphi(\overline{x})$, where $\overline{x}$ is an $m$-tuple of variables
containing all free variables in $\varphi$,
and any $\overline{a}\in\mathcal{M}^m$, called parameters, we can define
$\mathcal{M}\models\varphi(\overline{a})$,
meaning $\varphi(\overline{a})$ is true in $\mathcal{M}$, recursively by (a)
$\mathcal{M}\models a=a'$ iff $a$ and $a'$ are identical elements
in $\mathcal{M}$ and $\mathcal{M}\models P_i(\overline{a})$ iff $\overline{a}\in
P_i^{\mathcal{M}}$; (b) for any $\mathscr{L}$--formulas $\varphi$ and $\psi$
with all free variables $\overline{x}$ being substituted
by $\overline{a}$ in $\mathcal{M}$,
$\mathcal{M}\models \neg\varphi$ iff
$\mathcal{M}\not\models\varphi$ ($\neg$ means ``not''),
$\mathcal{M}\models\varphi\vee\psi$ iff
$\mathcal{M}\models\varphi$ or $\mathcal{M}\models\psi$ ($\vee$ means ``or''),
$\mathcal{M}\models\varphi\wedge\psi$ iff
$\mathcal{M}\models\varphi$ and $\mathcal{M}\models\psi$, ($\wedge$ means ``and''),
$\mathcal{M}\models\varphi\to\psi$ iff $\mathcal{M}\models\neg\varphi\vee\psi$
($\to$ means ``imply''),
and $\mathcal{M}\models\varphi\leftrightarrow\psi$
iff $\mathcal{M}\models(\varphi\to\psi)\wedge(\psi\to\varphi)$ ($\leftrightarrow$
means ``if and only if'');
(c) for any $\mathscr{L}$--formulas $\varphi(y,\overline{x})$,
$\mathcal{M}\models\forall y\,\varphi(y,\overline{a})$
iff $\mathcal{M}\models\varphi(b,\overline{a})$
for every $b\in\mathcal{M}$ ($\forall$ means
``for all'') and $\mathcal{M}\models\exists y\varphi(y,\overline{a})$ iff
$\mathcal{M}\models\varphi(b,\overline{a})$ for some $b\in\mathcal{M}$
($\exists$ means ``there exist''). From the truth definition above, for
every $\mathscr{L}$--formula $\varphi(\overline{x})$
there is another $\mathscr{L}$--formula $\psi(\overline{x})$
using only logic connectives $\neg$ and $\vee$, and only quantifier $\exists$,
such that $\mathcal{M}\models\varphi(\overline{a})$ iff
$\mathcal{M}\models\psi(\overline{a})$ for any $\mathscr{L}$--model
$\mathcal{M}$ and any $\overline{a}\in\mathcal{M}^m$.

We sometimes call a formula with all free variables being substituted by
parameters a sentence. Clearly, the truth value of a sentence with parameters
from $\mathcal{M}$ is determined in $\mathcal{M}$.
When we write $\varphi(\overline{x},\overline{a})$, we mean implicitly
that all free variables in the formula are among $\overline{x}$
and all parameters from model $\mathcal{M}$ are among $\overline{a}$.

Suppose $\mathcal{M}$ and $\mathcal{N}$ are two $\mathscr{L}$--models.
A function $i:\mathcal{M}\to\mathcal{N}$
is called an {\it elementary embedding} from $\mathcal{M}$
to $\mathcal{N}$ if for any $\mathscr{L}$--formula $\varphi(\overline{x})$
and any $m$-tuple $\overline{a}=(a_1,a_2,\ldots,a_m)\in\mathcal{M}^m$ we have
\begin{equation}\label{transfer}
\mathcal{M}\models\varphi(\overline{a})\,\mbox{ iff }\,
\mathcal{N}\models\varphi(i(\overline{a}))
\end{equation}
where $i(\overline{a}):=(i(a_1),i(a_2),\ldots,i(a_m))\in\mathcal{N}^m$.
An elementary embedding is necessarily injective.
If there exists an elementary embedding $i:\mathcal{M}\to\mathcal{N}$,
we can view $\mathcal{M}$ as an elementary submodel of $\mathcal{N}$
 and call $\mathcal{N}$ an elementary
extension of $\mathcal{M}$, denoted by $\mathcal{M}\preceq\mathcal{N}$.
We sometimes write $\mathcal{M}\prec\mathcal{N}$
to emphasize that $i$ is not surjective.

If $\mathcal{M}$ is an $\mathscr{L}$--model,
$\mathscr{L}':=\mathscr{L}\cup\{P_{n+1},\ldots,P_k\}$, and $P_i^{\mathcal{M}}
\subseteq\mathcal{M}^{m_i}$ for $n<i\leq k$, the $\mathscr{L}'$--model $\mathcal{M'}$
by adding the relations $P_i^{\mathcal{M}}$ to $\mathcal{M}$ is denoted by
$(\mathcal{M};P^{\mathcal{M}}_{n+1},\ldots,P^{\mathcal{M}}_k)$. We
call $\mathcal{M}'$ a {\it model expansion} of $\mathcal{M}$.

\bigskip

\noindent {\bf \S\ref{nsa}.3 Ultrapower construction}\quad
Next we construct an elementary extension of a model using ultrapower
construction (cf.\ \cite[\S4]{CK}).

\begin{definition}\label{ultrafilter}
Let $X$ be an infinite set. A set $\f\subseteq\mathscr{P}(X)$ is a
non-principal ultrafilter on $X$ if it satisfies the following:
\begin{enumerate}
\item $X\in\f$ and $F\not\in\f$ for any $F\in\mathscr{P}_{<\omega}(X)$
where $\mathscr{P}_{<\omega}(X)$ is the collection of all finite subsets of $X$,
\item $\forall A,B\in\mathscr{P}(X)\,(A,B\in\f\,\mbox{ implies }\,A\cap B\in\f)$,
\item $\forall A,B\in\mathscr{P}(X)\,(A\in\f\,\mbox{ and }\, A\subseteq B
\,\mbox{ imply }\, B\in\f)$,
\item $\forall A\in\mathscr{P}(X)\,(A\in\f\,\mbox{ or }\, (X\setminus A)\in\f)$.
\end{enumerate}
\end{definition}

It is well known that the existence of
non-principal ultrafilters on $X$ follows from $\mathsf{ZFC}$.

\begin{definition}
Let $\mathcal{M}$ be an $\mathscr{L}$--model and
$\f$ be a non-principal ultrafilter on an infinite set $X$.
Let $\mathcal{M}^X$ be the set of all functions from $X$ to $\mathcal{M}$.
For any $f,g\in\mathcal{M}^X$ define
\[f\sim_{\f} g\,\mbox{ iff }\,\{n\in X\mid f(n)=g(n)\}\in\f.\]
It is easy to check that $\sim_{\f}$ is an equivalence relation. Denote
$[f]_{\f}$ for the equivalence class containing $f$. The ultrapower
of $\mathcal{M}$ modulo $\f$, denoted by $\mathcal{M}^{X}/\f$,
is an $\mathscr{L}$--model which is composed of the base set
\[\{[f]_{\f}\mid f\in\mathcal{M}^{X}\}\] and the interpretation of $P_i$ by
\[P_i^{\mathcal{M}^X/\f}:=\left\{([f_1]_{\f},[f_2]_{\f},
\ldots,[f_{m_i}]_{\f})\mid
\{n\in X\mid (f_1(n),f_2(n),\ldots,f_{m_i}(n))
\in P_i^{\mathcal{M}}\}\in\f\right\}\]
for every $P_i\in\mathscr{L}$. Notice that
\[P_i^{\mathcal{M}^X/\f}=\{[\overline{f}]_{\f}\mid \overline{f}\,\mbox{ is a function
from }\,X\,\mbox{ to }\,P_i^{\mathcal{M}}\}.\]
The right side above is the ultrapower of $P_i^{\mathcal{M}}$ modulo $\f$.
\end{definition}

For an element $c\in\mathcal{M}$ denote $\phi_c:X\to\mathcal{M}$
for the constant function with value $c$.
Let $i:\mathcal{M}\to\mathcal{M}^X/\f$ be the {\it natural embedding},
i.e.,  $i(c)=[\phi_c]_{\f}$. The following is often called {\L}o\'{s}'s theorem.
(cf.\ \cite[Theorem 4.1.9, Corollary 4.1.13]{CK}.)

\begin{proposition}\label{los}
For any $\mathscr{L}$--formula $\varphi(\overline{x})$ and
any $\overline{[a]}_{\f}\in\mathcal{M}^X/\f$ it is true that
\[\mathcal{M}^X/\f\models\varphi(\overline{[a]}_{\f})\,\mbox{ iff }\,
\left\{n\in X\mid\mathcal{M}\models\varphi(\overline{a(n)})\right\}\in\f.\]
\end{proposition}

The proof of Proposition \ref{los} is done by induction on the complexity
of $\varphi$.

\begin{corollary}
The natural embedding $i:\mathcal{M}\to\mathcal{M}^X/\f$ with $i(c)=[\phi_c]_{\f}$ is
an elementary embedding. Furthermore, if $\mathcal{M}'=(\mathcal{M};R)$
is a model expansion of $\mathcal{M}$, then $i$ is also an elementary embedding
from $\mathcal{M}'$ to $\mathcal{M}'^X/\f=(\mathcal{M}^X/\f;R^X/\f)$.
\end{corollary}

So, the model $\mathcal{M}$ can be viewed
as an elementary submodel of $\mathcal{M}^X/\f$ via the natural embedding $i$
and $\mathcal{M}^X/\f$ is an elementary extension of $\mathcal{M}$.

\bigskip

\noindent {\bf \S\ref{nsa}.4 Construction of $\V_1$}\quad
Fix a non-principal ultrafilter $\f_0$ on $\n_0$.

\medskip

$\blacklozenge$: From now on let $\mathscr{L}=\{\in\}$. The ultrafilters
which will be used are $\f_0$, the tensor product of $\f_0$'s, and
nonstandard versions of them.

\medskip

Recall that $\V_0$ represents the standard universe,
which is an $\mathscr{L}$--model. Let
$\V_0^{\n_0}/\f_0$ be the ultrapower of $\V_0$ modulo $\f_0$
and $i^0:\V_0\to\V_0^{\n_0}/\f_0$ be the natural embedding.
Denote $^*\!\!\in$ for the interpretation of $\in$ in $\V_0^{\n_0}/\f_0$.
Notice that $a\in b$ iff $i^0(a)\,^*\!\!\in i^0(b)$ for any $a,b\in\V_0$.
One can define the rank function $\mbox{rank}(b)$ for every $b\in\V_0^{\n_0}/\f_0$
according to (\ref{rank}) with $\in$ being replaced by $^*\!\!\in$.
Notice that some elements $[f]_{\f}$ in $\V_0^{\n_0}/\f_0$ may not have
a finite $^*\!\!\in$-rank. For example, if $f(1)=0$ and $f(n+1)=\{f(n)\}$
for every $n\in\n_0$, then $[f]_{\f_0}\in\V_0^{\n_0}/\f_0$ does not have a
finite $^*\!\!\in$-rank. Let
\[V_1:=\left\{b\in\V_0^{\n_0}/\f_0\mid\,^*\!\!\in\!\mbox{-rank}(b)\in\omega\right\}.\]
The set $V_1$ is just the ultrapower of $\V_0$ modulo $\f_0$ truncated
at $^*\!\!\in$-rank $\omega$.
Notice that $\in\!\mbox{-rank}(a)=\,^*\!\!\in\!\mbox{-rank}(i^0(a))$
for every $a\in\V_0$ and
$V_1=\bigcup_{n\in\omega}i^0(V(\real_0,n))$.

Let $\real_1:=i^0(\real_0)$ and $\n_1:=i^0(\n_0)$.
Assume that every element in $\real_1$
is an urelement and identify $i^0(r)$ by $r$ for every $r\in\real_0$.
Then $\real_0\subseteq\real_1$.
Since the natural order $\leq$ on $\real_0$,
addition $+$ and multiplication $\times$ on $\real_0$
can be viewed as elements in $\V_0$,
we have that $i^0(\leq)$ is a linear order on $\real_1$ extending $\leq$,
$i^0(+)$ is the addition on $\real_1$ extending $+$, and $i^0(\times)$ is the
multiplication on $\real_1$ extending $\times$. For notational convenience
we write $\leq$, $+$, and $\times$ for $i^0(\leq)$, $i^0(+)$, and $i^0(\times)$,
respectively. By the elementality of $i^0$
the structure
$(\real_1;+,\times,\leq,0,1)$ is an ordered field containing the standard
real field $(\real_0;+,\times,\leq,0,1)$ as its subfield.
Notice that if $\mbox{Id}(n)=n$ for every $n\in\n_0$, then
$[\mbox{Id}]_{\f_0}\in i^0(\n_0)=\n_1$
and $[\mbox{Id}]_{\f_0}\geq r$ for every $r\in\real_0$, i.e., $\n_1$ contains
natural numbers such as $[\mbox{Id}]_{\f_0}$ which are infinitely large
relative to real numbers in $\real_0$.

Let $\mathscr{M}$ be the Mostowski collapsing map on $V_1$, i.e.,
$\mathscr{M}(a)=a$ for every $a\in \real_1$ and
\[\mathscr{M}(b):=\{\mathscr{M}(a)\mid a\,^*\!\!\in b\}\]
for every $b\in V_1\setminus \real_1$. Then $\mathscr{M}$
is an injection and $a\,^*\!\!\in b$ iff $\mathscr{M}(a)\in\mathscr{M}(b)$.
If one identifies $V_1$ with the image of $V_1$ under $\mathscr{M}$,
one can pretend that $^*\!\!\in$ is the true membership relation
and consider $V_1$ as a subset of the superstructure $\V(\real_1)$.
Hence, we can pretend that $^*\!\!\in$ is the true membership relation
$\in$ and drop $^*$ for notational convenience.
The purpose of the truncation of $\V_0^{\n_0}/\f_0$ at $^*\!\!\in$-rank $\omega$
is to make sure $\mathscr{M}$ is well defined in the standard sense.

Let $\V_1:=(V_1;\in)$. We call $\V_1$ a {\it nonstandard universe}
 extending $\V_0$. We will extend $\V_1$ further later. Notice that
 due to the truncation, $\V_1$ is no longer an elementary extension of
 $\V_0$ from the model theoretic point of view. However, $\V_1$
 is a so-called {\it bounded elementary extension} of $\V_0$.

An $\mathscr{L}$--formula
$\theta$ has {\it bounded quantifiers} if every occurrence of quantifiers $\forall$
and $\exists$ in $\theta$ has the form $\forall x\!\in\!y$
and $\exists x\!\in\!y$. Notice that $\forall x\!\in\!y\,\varphi$
is the abbreviation of $\forall x(x\!\in\!y\to\varphi)$ and
$\exists x\!\in\!y\,\varphi$ is the abbreviation of
$\exists x(x\!\in\!y\wedge\varphi)$. Similar to (\ref{transfer}),
it is easy to show that for any $\mathscr{L}$--formula
$\varphi(\overline{x})$ with bounded quantifiers
and any $\overline{a}\in\V_0$ we have
\begin{equation}\label{transfer2}
\V_0\models\varphi(\overline{a})\,\mbox{ iff }\,
\V_1\models\varphi(i^0(\overline{a})).
\end{equation}
So, the map $i^0$ is called a {\it bounded elementary embedding}
from $\V_0$ to $\V_1$. It is a common abuse of notation to write
$\V_0\preceq\V_1$
to indicate the existence of the bounded elementary embedding (instead of
just elementary embedding) $i^0$ and $\V_0\prec\V_1$ to
emphasize that $i^0$ is not surjective. The property (\ref{transfer2})
is sometimes called the {\it transfer principle}
in nonstandard analysis. Notice that if $A\in\V_0$, then
$i^0(A)$ can be viewed as the ultrapower of $A$ modulo $\f_0$, i.e.,
$i^0(A)=\{[f]_{\f_0}\mid f\in A^{\n_0}\}$.

Each $i^0(a)$ for $a\in\V_0$ is called $\V_0$--internal (or ``standard'' in
some literature) and
each $b\in\V_1$ is called $\V_1$--internal. Hence $\V_0$--internal set
is also a $\V_1$--internal set.
For example, $\real_1=i^0(\real_0)$ and
$\n_1=i^0(\n_0)$ are $\V_0$--internal sets.
Some $\V_1$--internal sets are not $\V_0$--internal. For example,
the set $\{1,2,\ldots,[\mbox{Id}]_{\f_0}\}$
is $\V_1$--internal subset of $\n_1$ but not $\V_0$--internal.
Some subsets of a $\V_1$--internal set are not
$\V_1$--internal. For example, $\n_0$ as a subset of $\n_1$ is
not $\V_1$--internal because it is bounded above in $\n_1$ and has no
largest element. Notice that $i^0[\n_0]=\n_0$ and $i^0(\n_0)=\n_1$.
The following proposition says that a subset of $B\in\V_1$
defined by an $\mathscr{L}$--formula with parameters
from $\V_1$ is $\V_1$--internal. The proposition is
an easy consequence of Proposition \ref{los}.

\medskip

\begin{proposition}\label{internal1}
Let $\varphi(\overline{a},\overline{x})$ be an $\mathscr{L}$--formula
with bounded quantifiers, and parameters $\overline{a}$ and $B^m$
being $\V_1$--internal. Then
$\left\{\overline{b}\in B^m\mid\V_1\models\varphi(\overline{a},\overline{b})\right\}$
is a $\V_1$--internal set. (cf.\ \cite[Theorem 4.4.14]{CK}.)
\end{proposition}

The following is called the {\it overspill principle}
in nonstandard analysis. Notice that $\n_0$ is an infinite
initial segment of $\n_1$.

\begin{proposition}\label{spill1}
Let $U\subseteq\n_1$ be an infinite proper
initial segment of $\n_1$ and not $\V_1$--internal.
Let $A$ be an $\V_1$--internal subset of $\n_1$.
If $A\cap U$ is upper unbounded in $U$,
then $A\cap(\n_1\setminus U)\not=\emptyset$.
\end{proposition}

The proof of Proposition \ref{spill1} is easy.
If $A\cap(\n_1\setminus U)=\emptyset$, then $U$
can be defined by a formula with bounded quantifiers and
parameter $A$. Hence $U$ is $\V_1$--internal.

\bigskip

\noindent {\bf \S\ref{nsa}.5 Construction of $\V_2$ and $\V_3$}\quad
We now extend $\V_1$ further to $\V_2$ and $\V_3$ to form a nonstandard
extension chain
\[\V_0\prec\V_1\prec\V_2\prec\V_3\] using ultrapower construction
and show the existence of bounded elementary embeddings $i_0$,
$i_*$, $i_1$, and $i_2$ besides the natural embeddings $i^1:\V_1\to\V_2$
and $i^2:\V_2\to\V_3$. The chain and embeddings
will satisfy the following properties. Let $\n_{j+1}=i^j(\n_j)$ and
$\real_{j+1}=i^j(\real_j)$ be the set of all positive
integers and the set of all real numbers, respectively, in $\V_{j+1}$ for $j=1,2$.

\begin{property}\label{property}
\begin{enumerate}
\item For $j=0,1,2$, $\n_{j+1}$ is an end--extension of $\n_j$, i.e.,
every number in $\n_{j+1}\setminus\n_j$ is greater than each number
in $\n_j$.

\item There is a bounded elementary embedding $i_*$ from the
$\mathscr{L}'$--model \newline $(\V_2;\real_0,\real_1)$ to
the $\mathscr{L}'$--model $(\V_3;\real_1,\real_2)$, where
$\mathscr{L}':=\mathscr{L}\cup\{P_1,P_2\}$ for
two new unary predicate symbols $P_1$ and $P_2$ not in $\mathscr{L}$.
Furthermore, the map $i_1:=i_*\res\V_1$ is a bounded elementary embedding
from $(\V_1;\real_0)$ to $(\V_2;\real_1)$. Notice that
$i_1(a)\in\n_2\setminus\n_1$ for each $a\in\n_1\setminus\n_0$.

\item There is a bounded elementary embedding $i_2$ from $\V_2$ to $\V_3$
such that $i_2\res\n_1$ is an identity map and $i_2(a)\in\n_3\setminus\n_2$
for each $a\in\n_2\setminus\n_1$.
\end{enumerate}
\end{property}

We are now going to work towards establishing this property in this section.

\medskip

$\blacklozenge$: From now on, let $\V_j^{\n_j}$
always represent, for notational convenience,
the set of all functions $f$ from $\n_j$ to $\V_j$ such that
$\{\mbox{rank}(f(n))\mid n\in\n_j\}$ is a {\bf bounded set} in $\omega$.

\medskip

Notice that $\f_1:=i^0(\f_0)\in\V_1$
satisfies Part 1--4 of Definition \ref{ultrafilter}
with $X$, $\mathscr{P}_{<\omega}(X)$, and $\mathscr{P}(X)$
being replaced by $\n_1:=i^0(\n_0)$,
$i^0(\mathscr{P}_{<\n_0}(\n_0))=\left\{A\subseteq\n_1\mid
A\in\V_1\,\mbox{ and }\,|A|\in\n_1\right\}$, and
$i^0(\mathscr{P}(\n_0))=\V_1\cap\mathscr{P}(\n_1)$, respectively.
We call $\f_1$ a $\V_1$--internal
non-principal ultrafilter on $\n_1$.

\begin{definition}\label{v2}
Let $\f_1:=i^0(\f_0)$. Denote $\V_2$ for the model $(V_2;\,\in_2)$ such that
\[V_2:=(\V_1^{\n_1}\cap\V_1)/\f_1\,\mbox{ and}\]
\[[f]_{\f_1}\in_2 [g]_{\f_1}\,\mbox{ iff }\,\{n\in\n_1\mid
f(n)\in g(n)\}\in\f_1\] for all $f,g\in\V_1^{\n_1}\cap\V_1$. Let
$i^1:\V_1\to\V_2$ be the natural embedding that $i^1(c)=[\phi_c]_{\f_1}$ for every
$c\in\V_1$.
\end{definition}

We call $\V_2$ a $\V_1$--internal ultrapower of $\V_1$ modulo
the $\V_1$--internal ultrafilter $\f_1$.

\begin{proposition}\label{los2}
If $\varphi(\overline{x})$ is an $\mathscr{L}$--formula with
bounded quantifier and \newline
$\overline{[a]}_{\f_1}\in\V_2$, then
\[\V_2\models\varphi(\overline{[a]}_{\f_1})\,\mbox{ iff }\,
\left\{n\in \n_1\mid\V_1\models\varphi(\overline{a(n)})\right\}\in\f_1.\]
\end{proposition}

\begin{corollary}
The natural embedding $i^1:\V_1\to\V_2$ is a bounded
elementary embedding from $\V_1$ to $\V_2$.
\end{corollary}

The proof of Proposition \ref{los2} is almost the same as the
proof of Proposition \ref{los} except one step that shows
$\left\{n\in\n_1\mid\V_1\models\exists x\!\in\!B(n)\,\varphi(\overline{a(n)},x)
\right\}\in\f_1$ implies $\V_2\models\exists x\!\in\! [B]_{\f_1}\,
\varphi(\overline{[a]}_{\f_1},x)$. Let $M_z=i^0(\V(\real_0,z))\in\V_1$ with
$z\in\omega$ such that $B,\overline{a}\in M_z$.
By the axiom of choice there is a well-order $\lhd$ on $\V(\real_0,z)$. So,
every nonempty set $A\in M_z$ has a $i^0(\lhd)$-least element
by the transfer principle. Suppose that
$\left\{n\in\n_1\mid\V_1\models\exists x\!\in\!B(n)\,\varphi(\overline{a(n)},x)
\right\}=X\in\f_1$. For each $n\in\n_1$, if $n\not\in X$, let $f(n)=0$ and if
$n\in X$, let $f(n)$ be the $i^0(\lhd)$-least element in
the nonempty set $\left\{x\in B(n)\mid\V_1\models\varphi(a(n),x)\right\}$.
Since $B,\overline{a},i^0(\lhd)$ are in $\V_1$, so does the function $f$
by Proposition \ref{internal1}.
Hence $[f]_{\f_1}\in\V_2\cap [B]_{\f_1}$. By the induction hypothesis we have that
\[\left\{n\in\n_1\mid\V_1\models\varphi(a(n),f(n))\right\}\supseteq X\in\f_1\,\mbox{
implies }\,\V_2\models\varphi(\overline{[a]}_{\f_1},[f]_{\f_1}),\] which implies
$\V_2\models\exists x\!\in\! [B]_{\f_1}\,
\varphi(\overline{[a]}_{\f_1},x)$.

Let
\begin{equation}\label{i_0}
i_0=i^1\!\circ\!i^0.
\end{equation}
Then $i_0$ is a bounded elementary embedding from $\V_0$ to $\V_2$,
which will be used in Theorem \ref{szemeredi}.

Same as for $\V_1$ we can assume by Mostowski collapsing that $\V_2$
is a subset of the superstructure $\V(\real_2)$ and
$\in_2$ is the true membership relation $\in$.
Notice that $i^1\res\real_1$ is an identity map.
The element $a\in\V_2$ is called $\V_2$--internal,
$i^1(b)$ is called $\V_1$--internal for any $b\in\V_1$, and $i_0(c)$
is called $\V_0$--internal for every $c\in\V_0$.
Notice that $\n_0$ and $\n_1$ as subsets of $\n_2$
are not $\V_2$--internal.

If $f\in\V_1$ is a
function from $\n_1$ to $[n]$ for some $n\in\n_1$, then
$[f]_{\f_1}=m$ for some $m\in\n_1$ by (\ref{transfer2}) with
$\V_0,\V_1,i^0$ being replaced by $\V_1,\V_2,i^1$.
Hence $\n_2$ is a proper end-extension of $\n_1$.

\medskip

By the same way, we can define $\V_3$ as a $\V_2$--internal ultrapower
of $\V_2$ modulo $i^1(\f_1)$.

\begin{definition}
Let $\f_2:=i^1(\f_1)$ be the $\V_2$--internal ultrafilter on $\n_2$.
Let $\V_3$ be the model $(V_3;\,\in_3)$ such that
\[V_3:=(\V_2^{\n_2}\cap\V_2)/\f_2\,\mbox{ and}\]
\[[f]_{\f_2}\in_3 [g]_{\f_2}\,\mbox{ iff }\,\{n\in\n_2\mid
f(n)\in g(n)\}\in\f_2\] for all $f,g\in\V_2^{\n_2}\cap\V_2$, and define
the natural mebedding $i^2:\V_2\to\V_3$ by $i^2(c)=[\phi_c]_{\f_2}$
for every $c\in\V_2$.
\end{definition}

Generalizing the arguments above, we have that the map $i^2$ is a bounded
elementary embedding from $\V_2$ to $\V_3$. We say that $\V_3$ is a
nonstandard extension of $\V_2$. It is also easy to see that
$\n_3:=i^2(\n_2)$ is an end-extension of $\n_2$.

We have completed the construction of $\V_0\prec\V_1\prec\V_2\prec\V_3$
and verified that Part 1 of Property \ref{property} is true.

\bigskip

\noindent {\bf \S\ref{nsa}.6 Bounded elementary embeddings $i_*$, $i_1$, and $i_2$}\quad
To verify Part 2--3 of Property \ref{property} we have to view the construction
of $\V_j$ for $j=2,3$ from a different angle. For a set $A\subseteq \n_0\times\n_0$
and $n\in\n_0$, let $A_n:=\left\{m\in\n_0\mid (m,n)\in A\right\}$.
Let $\f_0$ and $\f_0'$ be two non-principal ultrafilters on $\n_0$.
The tensor product of $\f_0$ and $\f'_0$ is defined by
\[\f_0\otimes\f_0':=\left\{A\subseteq\n_0
\times\n_0\mid \{m\in\n_0\mid A_m\in\f_0\}\in\f_0'\right\}.\]
It is easy to check that $\f_0\otimes\f_0'$ is a non-principal
ultrafilter on $\n_0\times\n_0$. For simplicity we assume that
$\f_0$ and $\f'_0$ are the same ultrafilter.

\begin{lemma}\label{v21}
Let $\V'_2:=(V_0^{\n_0\times\n_0}/\f_0\otimes\f_0';\in')$ where
$[f]_{\f_0\otimes\f'_0}\in' [g]_{\f_0\otimes\f'_0}$ iff\newline
$\left\{(m,n)\in\n_0\times\n_0\mid f(m,n)\in g(m,n)\right\}\in\f_0\otimes\f_0'$. Then
\[\V_2\cong\V_2'.\]
\end{lemma}

\noindent {\it Proof}\quad Let $m$ range over the first copy of
$\n_0$ and $n$ over the second copy of $\n_0$ in $\n_0\times\n_0$.
Given $a\in\V_2$, there is an $f_a\in\V_1^{\n_1}\cap\V_1$
such that $a=[f_a]_{\f_1}$. Since the ranks of all image of $f_a$ is bounded, the range
of $f_a$ is in a $\V_0$--internal set $i^0(B)$.
So, by identifying $f_a$ with its graph, we have $f_a\subseteq\n_1\times i^0(B)$.
Since $f_a\in\V_1$, there is a $g_a\in\V_0^{\n_0}$ such that
$f_a=[g_a]_{\f_0}$ where $g_a(n)\subseteq\n_0\times B$
is the graph of a function $g_a(n):\n_0\to B$ for all $n\in\n_0$.
Now let $F_a:\n_0\times\n_0\to B$ be such that
$F_a(m,n)=g_a(n)(m)\in B$. Then $F_a\in\V_0^{\n_0\times\n_0}$.

We view $a\mapsto [F_a]_{\f_0\otimes\f'_0}$ as a relation between $\V_2$
and $\V'_2$. Notice that $[X]_{\f'_0}\in\f_1$ iff
$\{n\in\n_0\mid X(n)\in\f_0\}\in\f'_0$.
Notice also that for any $a,b\in\V_2$, we have
\begin{eqnarray*}
\lefteqn{a=b\,\mbox{ iff }}\\
& &\left\{[x]_{\f'_0}\in\n_1\mid f_a([x]_{\f'_0})
=f_b([x]_{\f'_0})\right\}\in\f_1\,\mbox{ iff}\\
& &\left\{[x]_{\f'_0}\in\n_1\mid [g_a]_{\f'_0}([x]_{\f'_0})=
[g_b]_{\f'_0}([x]_{\f'_0})\right\}\in\f_1\,\mbox{ iff}\\
& &\left\{[x]_{\f'_0}\in\n_1\mid [g_a(n)(x(n))]_{\f'_0}=
[g_b(n)(x(n))]_{\f'_0}\right\}\in\f_1\,\mbox{ iff}\\
& &\left\{n\in\n_0\mid \{m=x(n)\in\n_0\mid g_a(n)(x(n))=
g_b(n)(x(n))\}\in\f_0\right\}\in\f'_0\,\mbox{ iff}\\
& &\left\{n\in\n_0\mid \{m\in\n_0\mid g_a(n)(m)=
g_b(n)(m)\}\in\f_0\right\}\in\f'_0\,\mbox{ iff}\\
& &\left\{n\in\n_0\mid \{m\in\n_0\mid F_a(m,n)=
F_b(m,n))\}\in\f_0\right\}\in\f'_0\,\mbox{ iff}\\
& &\left\{(m,n)\in\n_0\times\n_0\mid F_a(m,n)=
F_b(m,n))\right\}\in\f_0\otimes\f'_0\,\mbox{ iff}\\
& &[F_a]_{\f_0\otimes\f'_0}=[F_b]_{\f_0\otimes\f'_0}.
\end{eqnarray*}
Hence, the relation
$a\mapsto [F_a]_{\f_0\otimes\f'_0}$ is an injective function from $\V_2$ to $\V'_2$.

On the other hand, given
$[F]_{\f_0\otimes\f'_0}\in\V_0^{\n_0\times\n_0}/\f_0\otimes\f'_0$, let $g(n)\in
\V_0^{\n_0}$ be such that $g(n)(m)=F(m,n)$.
Let $f:\n_1\to\V_1$ be such that $f([x]_{\f_0}):
=[g]_{\f'_0}([x]_{\f_0})=[g(n)(x(n))]_{\f'_0}$. Then $f=[g]_{\f'_0}\in\V_1$.
So, $a=[f]_{\f_1}\in\V_2$ is such that
$[F_a]_{\f_0\otimes\f'_0}=[F]_{\f_0\otimes\f'_0}$. This shows that
$a\mapsto [F_a]_{\f_0\otimes\f'_0}$ is surjective.

Notice also that for any $a,b\in\V_2$,
\begin{eqnarray*}
\lefteqn{a\in b\,\mbox{ iff}}\\
& &\left\{[x]_{\f'_0}\in\n_1\mid f_a([x]_{\f'_0})\in
f_b([x]_{\f'_0})\right\}\in\f_1\,\mbox{ iff}\\
& &\left\{[x]_{\f'_0}\in\n_1\mid [g_a]_{\f'_0}([x]_{\f'_0})\in
[g_b]_{\f'_0}([x]_{\f'_0})\right\}\in\f_1\,\mbox{ iff}\\
& &\left\{[x]_{\f'_0}\in\n_1\mid [g_a(n)(x(n))]_{\f'_0}\in
[g_b(n)(x(n))]_{\f'_0}\right\}\in\f_1\,\mbox{ iff}\\
& &\left\{n\in\n_0\mid \{m=x(n)\in\n_0\mid g_a(n)(x(n))\in
g_b(n)(x(n))\}\in\f_0\right\}\in\f'_0\,\mbox{ iff}\\
& &\left\{n\in\n_0\mid \{m\in\n_0\mid g_a(n)(m)\in
g_b(n)(m)\}\in\f_0\right\}\in\f'_0\,\mbox{ iff}\\
& &\left\{n\in\n_0\mid \{m\in\n_0\mid F_a(m,n)\in
F_b(m,n))\}\in\f_0\right\}\in\f'_0\,\mbox{ iff}\\
& &\left\{(m,n)\in\n_0\times\n_0\mid F_a(m,n)\in
F_b(m,n))\right\}\in\f_0\otimes\f'_0\,\mbox{ iff}\\
& &[F_a]_{\f_0\otimes\f'_0}\in' [F_b]_{\f_0\otimes\f'_0}.
\end{eqnarray*}
This completes the proof of $\V_2\cong\V'_2$.
\hfill $\blacksquare$

\medskip

If we identify each of $\V_2$ and $\V'_2$
with its image of Mostowski collapsing, then  $\V_2$ and $\V'_2$ can be viewed
as the same model.

\begin{lemma}\label{v22}
Let $\V''_2:=\left((V_0^{\n_0}/\f_0)^{\n_0}/\f_0';\in''\right)$ where
$[[f]_{\f_0}]_{\f'_0}\in''[[g]_{\f_0}]_{\f'_0}$ iff\newline
$\left\{n\in\n_0\mid \left\{m\in\n_0\mid f(m,n)\in g(m,n)\right\}\in\f_0\right\}
\in\f_0'$. Then \[\V'_2\cong\V_2''\] for any $f,g:\n_0\times\n_0\to\V_0$.
\end{lemma}

The proof of Lemma \ref{v22} can be found in \cite[Proposition 6.5.2]{CK}.
We call $\V_2''$ the external ultrapower of $\V_1$ modulo $\f'_0$. By Lemma \ref{v21}
and Lemma \ref{v22} we can view $\V_2$ and $\V_2''$ as the same model and write
$\in$ for $\in''$. To summarize, we have that
\begin{equation}\label{twosteps}
\V''_2=(\V_0^{\n_0}/\f_0)^{\n_0}/\f_0'=\V_0^{\n_0\times\n_0}/\f_0\otimes\f_0'
=(\V_1^{\n_1}\cap\V_1)/\f_1=\V_2.
\end{equation}
The term on the left side is the two-step iteration when we take the ultrapowers by using
$\f_0$ first and $\f'_0$ second. The term on the right side is when we use
$\f'_0$ first so that $V_0^{\n_0}$ becomes $\V_1$, $\n_0$ becomes $\n_1$,
and $\f_0$ becomes $\f_1$.

However, the bounded elementary embedding from $\V_1$ to
$\V_2$ induced by the $\V_1$--internal ultrapower modulo $\f_1$ and the
bounded elementary embedding from $\V_1$ to $\V_2''$
induced by the external ultrapower of $\V_1$ modulo $\f_0'$ are different.
Let $i_1:\V_1\to\V_2''=\V_1^{\n_0}/\f'_0$ be such that $i_1(c)=[\phi_c]_{\f'_0}$
for each $c\in\V_1$.
If $c\in\n_0\subseteq\n_1$, then clearly we have $[\phi_c]_{\f_0'}=c$.
If $c,c'\in\n_1\setminus\n_0$ with $c=[f]_{\f'_0}$ and $c'=[g]_{\f'_0}$
for some $f,g\in\n_0^{\n_0}$, then $c>g(n)$ for every $n\in\n_0$.
Hence $i_1(c)=[\phi_c]_{\f_0'}>[g]_{\f_0'}=c'$.
This shows that $i_1(c)\in\n_2\setminus\n_1$. The map $i_1$ will
be $i_*\res\V_1$ for a bounded elementary embedding $i_*$ defined below.
Notice that above arguments still work if the ultrafilters
$\f_0$ and $\f_0'$ are different.

\medskip

Generalizing the construction further and with the help of Mostowski
collapsing, we have that
\begin{eqnarray}
\lefteqn{\V''_3:=((\V_0^{\n_0}/\f_0)^{\n_0}/\f_0')^{\n_0}/\f_0''
=\V_2^{\n_0}/\f''_0}\label{v31}\\
& &=\V_0^{\n_0\times\n_0\times\n_0}/\f_0\otimes\f_0'\otimes\f_0''\nonumber\\
& &=(\V_1^{\n_1\times\n_1}\cap\V_1)/\f_1\otimes\f'_1
=(\V_2^{\n_1}\cap\V_1)/\f'_1=\V_3'\label{v32}\\
& &=(\V_1^{\n_1\times\n_1}\cap\V_1)/\f_1\otimes\f'_1=
(\V_2^{\n_2}\cap\V_2)/\f_2=\V_3.\nonumber
\end{eqnarray}

Let $i_*:\V_2\to\V_3''=\V_2^{\n_0}/\f_0''$ be the bounded elementary embedding
determined by (\ref{v31}) with $i_*(c)=[\phi_c]_{\f_0''}$ for every
$c\in\V_2$. Since $\real_0^{\n_0}/\f''_0=\real_1$ and $\real_1^{\n_0}/\f''_0
=(\real_1^{\n_1}\cap\V_1)/\f_1=\real_2$ we conclude that the map
$i_*$ is also a bounded elementary embedding from
model expansion $(\V_2;\real_0,\real_1)$ to $(\V_3;\real_1,\real_2)$.
Since $\V_2=\V_1^{\n_0}/\f_0''$, we have that $i_*\res\V_1$ coincides with
$i_1$ mentioned above when $\V_2''$ is constructed, and
is a bounded elementary embedding from $(\V_1;\real_0)$ to $(\V_2;\real_1)$.
Hence Part 2 of Property \ref{property} is verified.

Let $i_2:\V_2\to\V_3'=(\V_2^{\n_1}\cap\V_1)/\f_1'$ be the bounded elementary embedding
determined by (\ref{v32}) with $i_2(c)=[\phi_c]_{\f_1'}$ for every
$c\in\V_2$. Since every $\V_1$--internal function from $\n_1$ to $[n]$ for some
$n\in\n_1$ is equivalent to a constant function modulo $\f'_1$ we conclude
that $i_2\res\n_1$ is an identity map. Similar to the argument for $i_1$, we have
$i_2(c)\in\n_3\setminus\n_2$ for every $c\in\n_2\setminus\n_1$.
This verifies Part 3 of Property \ref{property}.

\medskip

The following propositions are $\V_j$--versions of Proposition \ref{internal1}
and Proposition \ref{spill1}. Let $0<j\leq 3$.

\begin{proposition}\label{internal2}
Let $\varphi(\overline{a},\overline{x})$ be an $\mathscr{L}$--formula with
bounded quantifiers, and parameters $\overline{a}$ and $B^m$ be $\V_j$--internal. Then
$\left\{\overline{b}\in B^m\mid\V_j\models\varphi(\overline{a},\overline{b})\right\}$
is a $\V_j$--internal set.
\end{proposition}

\begin{proposition}\label{spill2}
Let $U$ be an infinite proper
initial segment of $\n_{j}$ and not $\V_j$--internal.
Let $A\subseteq\n_j$ be $\V_j$--internal.
If $A\cap U$ is upper unbounded in $U$,
then $A\cap(\n_j\setminus U)\not=\emptyset$.
\end{proposition}

We would like to mention that $\V_1$, $\V_2$, and $\V_3$ are ultrapowers of $\V_0$ modulo
the ultrafilters $\f_0$, $\f_0\otimes\f_0'$, and $\f_0\otimes\f'_0\otimes\f''_0$,
respectively. Hence, they are all countably saturated (cf.\ \cite[Corollary 4.4.24]{CK})
although the countable saturation is not used in this paper.
The proofs of Proposition \ref{internal2} and Proposition \ref{spill2}
are the same as Proposition \ref{internal1} and Proposition \ref{spill1}.
Notice also that $\V_2$ and $\V_3$ are $\V_1$--internal ultrapowers of
$\V_1$ modulo the $\V_1$--internal ultrafilters $\f_1$ and $\f_1\otimes\f'_1$,
respectively.

The ultrafilters $\f_0$, $\f'_0$, and $\f''_0$ do not have to be on the
countable set $\n_0$ and do not have to be the same. The only restriction
is that they have to be in $\V_0$.

Iterated nonstandard extensions were
used in combinatorial number theory before, e.g.\ in \cite{dinasso,DGL,luperi}.
But we will use them in a new way by exploring the advantages of
various bounded elementary embeddings between $\V_j$ and $\V_{j'}$ (see
Property \ref{property}). For most of the time we
work within $\V_2$. The nonstandard universe $\V_3$ is only used for
one step in the proof of Lemma \ref{key}.

\section{Nonstandard Versions of Some Facts}\label{nsatools}

\quad The Greek letters $\alpha,\beta,\eta,\epsilon$ will represent
``{\it standard}'' reals unless specified otherwise. Let $0\leq j<j'\leq 3$
in this paragraph.
An $r\in\real_3$ is a $\V_j$--infinitesimal, denoted by
$r\approx_j 0$, if $|r|<1/n$ for every $n\in\n_j$.
By an infinitesimal we mean a $\V_0$--infinitesimal. Denote $st_j$
for the $\V_j$--standard part map, i.e., $st_j(r)$ is the unique real number
$r'\in\real_j$ such that $r-r'\approx_j 0$ when $r\in\real_3\cap (-m,m)$
for some $m\in\n_j$. Notice that $st_j$ and $\n_j$ are definable
by a formula with bounded quantifiers and parameters in $(\V_{j'};\real_j)$.
Sometimes, the subscript $0$ will be dropped. For example,
$\approx$ means $\approx_0$ and $st$ means $st_0$.
For any two positive integers $m,n\in\n_3$ we denote $m\ll n$ for
$m\in\n_{j}$ and $n\in\n_{j'}\setminus\n_{j}$.
Hence $m\gg 1$ means $m\in\n_3\setminus\n_0$.
Notice that if $r\in\real_1\cap (-m,m)$ for some $m\in\n_0$ and $st(r)=\alpha$,
then $st(i_1(r))=\alpha$ where $i_1$ is in Part 2 of Property \ref{property}.
This is true because $i_1(\alpha)=\alpha$,
$i_1(1/n)=1/n$, and $|r-\alpha|<1/n$ iff
$|i_1(r)-i_1(\alpha)|<i_1(1/n)$ for each $n\in\n_0$.
Similarly, we have $st_1(i_2(r))=st_1(r)$
for $r\in\real_2$ in the domain of $st_1$.

Capital letters $A$, $B$, $C$, $\ldots$ represent sets of integers
except $H$, $J$, $K$, $N$ which are
reserved for integers in $\n_3\setminus\n_0$.
The letter $k\geq 3$ represents exclusively the
length of the arithmetic progression in
Szemer\'{e}di's theorem and $l$ represents an integer between $1$ and $k$.
All unspecified sets mentioned in this paper will be either standard subsets of
$\n_0$ or $\V_j$--internal sets for some $j=1,2,3$.
For any $n\in\n_3$ let $[n]:=\{1,2,\ldots,n\}$.

For any bounded set $A\subseteq\n_{j'}$ and $n\in\n_{j'}$
denote $\delta_n(A)$ for the quantity $|A|/n$ in $\V_{j'}$
where $|A|$ means the internal cardinality of $A$ in $\V_{j'}$. Denote
$\mu^j_n(A):=st_j(\delta_n(A))$ for $j<j'$.
Notice that $\delta_n$ is an internal function while $\mu^j_n$
are often external functions but definable in $(\V_{j'};\real_j)$
for $j'>j$, i.e.,
\[\mu^j_n(A)=\alpha\,\mbox{ iff }\,\forall n\in\n_{j'}\cap\real_j\,
\left(|\delta_H(A)-\alpha|<\frac{1}{n}\right).\] We often write
$\mu_n$ for $\mu^0_n$. If $A\subseteq\Omega$ and $|\Omega|=H$,
 then $\mu_H(A)$ coincides with the Loeb measure of
$A$ in $\Omega$. The term $\delta_H$ is often used for
an internal argument.

\medskip

$\blacklozenge$: The abbreviation {\it a.p.} stands for ``arithmetic progression''
and {\it $n$--a.p.} stands for ``$n$-term arithmetic progression.''

\medskip

The length of an a.p.\ $p$ is the number of terms in $p$ which
can be written as $|p|$. The letters $P,Q,R$ are reserved exclusively for
a.p.'s of length $\gg 1$, and $p,q,r$ for a.p.'s of
length $k$ or other standard length. When we run out of letters,
we may also use $\vec{x},\vec{y}$ for $k$--a.p.'s. If $1\leq l\leq |p|$, then
$p(l)$ represents the $l$-th term of $p$.
We denote $p\subseteq A$ for $p(l)\in A$ for all $1\leq l\leq |p|$.
We allow the common difference $d$ of an a.p.\ to be any
integer including, occasionally, the trivial case for $d=0$.
If $p$ and $q$ are two a.p.'s of the same length, then
$p\oplus q$ represents the $|p|$--a.p.\
$\{p(l)+q(l)\mid l=1,2,\ldots,|p|\}$. If $p$ is an a.p.\ and $X$
is an element or a set, then $p\oplus X$ represents the
sequence $\{p(l)+X\mid 1\leq l\leq |p|\}$. By $p\sqsubseteq q\oplus X$
we mean $p(l)\in q(l)+X$ for every $1\leq l\leq |p|=|q|$.

If $X\subseteq\real_j$ and $X\in\V_j$ let ${\sup}_j(X)$ be
the supremum of $X$ in the sense of $\V_j$, i.e., the unique least upper
bound of $X$ in $\real_j$, or $\infty$ if $X$ is unbounded
above in $\real_j$.

Let $1\leq j<j'\leq 3$. For any $A\subseteq\n_{j'-j}$ and
any collection of a.p.'s $\p\in\V_{j'-j}$,
there exists $X\in\V_0$ with $X\subseteq\real_0$
(because every subset of $\real_0$ is
in $\V_0$) such that $x\in X$ iff there exists a
$P\in\p$ with $|P|\in\n_{j'-j}\setminus\n_0$ and $\mu_{|P|}(A\cap P)=x$.
By the elementality of $i_*$ in Part 2 of Property \ref{property} we have that
for any $A\subseteq\n_{j'}$ and
any collection of a.p.'s $\p\in\V_{j'}$,
there exists $X\in\V_j$ with $X\subseteq\real_j$ such that $x\in X$ iff there exists a
$P\in\p$ with $|P|\in\n_{j'}\setminus\n_j$ and $\mu^j_{|P|}(A\cap P)=x$.
Therefore, the operator ${\sup}_j$ and hence $\SD^j$ below are well defined.
Similarly, $\SD_S^j$ below is also well defined.

\medskip

\begin{definition}\label{density}
For $0\leq j<j'\leq 3$ and $A\subseteq\n_{j'}$ with $|A|\in\n_{j'}\setminus
\n_j$ the strong upper
Banach density $\SD^j(A)$ of $A$ in $\V_j$ is defined by
\begin{equation}\label{sd-1}
\SD^j(A):={\sup}_j\left\{\mu^j_{|P|}(A\cap P)\mid |P|\in\n_{j'}\setminus\n_j\right\}.
\end{equation}
The letter $P$ above always represents an a.p.
If $S\subseteq\n_{j'}$ has $\SD^j(S)=\eta\in\real_j$ and $A\subseteq\n_{j'}$, the
strong upper Banach density $\SD_S^j$ of $A$ relative to $S$
is defined by
\begin{equation}\label{sds-1}
\SD^j_{S}(A):={\sup}_j\left\{\mu^j_{|P|}(A\cap P)\mid
|P|\in\n_{j'}\setminus\n_j,\mbox{ and }\mu^j_{|P|}(S\cap P)=\eta\right\}.
\end{equation}
\end{definition}

When $\SD^j_{S}(A)$, defined in (\ref{sds-1}), is used in this paper,
the set $A$ is often a subset of $S$ although there is no such restriction
in the definition.

\begin{definition}
If $A\subseteq\n_0$, then the strong upper Banach density
of $A$ is defined by $\SD(A):=\SD^0(i_0(A))$ where $\SD^0$ is defined
by (\ref{sd-1}) and $i_0$ is defined by (\ref{i_0}).
\end{definition}

We would like to point out that for standard sets $A\subseteq\n_0$,
\[\SD(A)=\lim_{n\in\n_0,\,n\to\infty}\sup\{\delta_{|p|}(A\cap p)\mid |p|\geq n\}.\]
This equality will not be used. The purpose here is to give the reader
some intuition because the right side is a standard expression.

The superscript $0$ in $\SD^0$ will be omitted.
Notice that the upper density of a set $A\subseteq\n_0$ is
less than or equal to the upper Banach density of $A$, which is less than
or equal to the strong upper Banach density of $A$.
The strong upper Banach density of $A$ is the nonstandard version of
the density of $A$ along a collection of arbitrarily long
arithmetic progressions satisfying the double counting property in \cite{tao}.
Of course, if we know that Szemer\'{e}di's theorem is true, then $\SD(S)=\eta>0$
implies that $\eta=1$. Also $\SD(S)=1$ and $\SD_S(A)=\alpha>0$ imply that $\alpha=1$.

Suppose that the strong upper Banach density of $A\subseteq\n_0$ is a positive
real number $\alpha$. Instead of looking for $k$--a.p.'s in $A$ we will look for
$k$--a.p.'s in $i_0(A)\cap P$ for some infinitely long a.p.\ $P$ such that
the distribution of $i_0(A)\cap P$ in $P$ is very uniform, i.e., the measure
and strong upper Banach density of $i_0(A)\cap P$ in $P$ are the same value $\alpha$.
The uniformity allows the use of an argument similar to
the so called density increment argument in the standard literature.
The next lemma is the beginning of this effort.

\begin{lemma}\label{equiv}
For $0\leq j<j'\leq 3$ let $A\subseteq S\subseteq\n_{j'}$
with $|A|\in\n_{j'}\setminus\n_j$
and $\alpha,\eta\in\real_j$ with $0\leq\alpha\leq\eta\leq 1$.
Then the following are true:
\begin{enumerate}
\item $\SD^j(S)\geq\eta$ iff there exists a $P$
with $|P|\in\n_{j'}\setminus\n_j$ and $\mu^j_{|P|}(S\cap P)\geq\eta$;
\item If $\SD^j(S)=\eta$, then there exists a $P$ with
$|P|\in\n_{j'}\setminus\n_j$ such that $\mu^j_{|P|}(S\cap P)=\SD^j(S\cap P)=\eta$;
\item Suppose $\SD^j(S)=\eta$. Then $\SD^j_S(A)\geq\alpha$ iff
 there exists a $P$ with $|P|\in\n_{j'}\setminus\n_j$,
$\mu^j_{|P|}(S\cap P)=\eta$, and $\mu^j_{|P|}(A\cap P)\geq\alpha$;
\item Suppose $\SD^j(S)=\eta$. If $\SD^j_S(A)=\alpha$, then
there exists a $P$ with $|P|\in\n_{j'}\setminus\n_j$
such that $\mu^j_{|P|}(S\cap P)=\eta$ and
$\mu^j_{|P|}(A\cap P)=\SD^j_{S\cap P}(A\cap P)=\alpha$.
\end{enumerate}
\end{lemma}

\noindent {\it Proof}\quad Part 1: If $\SD^j(S)\geq\eta$, then there is a $P_n$ with
$|P_n|\in\n_{j'}\setminus\n_j$ such that $\delta_{|P_n|}(S\cap P_n)>\eta-1/n$ for
every $n\in\n_j$. Let
\[A:=\left\{n\in\n_{j'}\mid \exists P\subseteq\n_{j'}\,(|P|\geq n\wedge
\delta_{|P|}(S\cap P)>\eta-1/n)\right\}.\] Then $A$ is $\V_{j'}$--internal
and $A\cap\n_j$ is unbounded above in $\n_j$. By Proposition \ref{spill2},
there is a $J\in A\setminus\n_j$. Hence there is an a.p.\
 $P_J\subseteq\n_{j'}$ such that $|P_J|\geq J\in\n_{j'}\setminus\n_j$
 and $\delta_{|P_J|}(S\cap P_J)>\eta-1/J\approx_j\eta$.
 Therefore, $\mu^j_{|P_J|}(S\cap P_J)\geq\eta$.
On the other hand, if $\mu^j_{|P|}(S\cap P)\geq\eta$, then
$\SD^j(S)\geq\eta$ by the definition of $\SD^j$ in (\ref{sd-1}).

\medskip

Part 2: If $\SD^j(S)=\eta$, we can find $P$ with $|P|\in\n_{j'}\setminus\n_j$
such that $\mu^j_{|P|}(S\cap P)=\eta'\geq\eta$ by Part 1. Clearly,
$\eta=\SD^j(S)\geq\SD^j(S\cap P)\geq\mu^j_{|P|}(S\cap P)
=\eta'$ by the definition of $\SD^j$. Hence $\eta=\eta'$.

\medskip

Part 3: If $\SD^j_S(A)\geq\alpha$, then there is a $P$ with $|P|>n$
such that $|\delta_{|P|}(S\cap P)-\eta|<1/n$ and
$\delta_{|P|}(A\cap P)>\alpha-1/n$ for every $n\in\n_j$.
By Proposition \ref{spill2} as in the proof of Part 1
there is a $P_J$ for some
$J\in\n_{j'}\setminus\n_j$ with $|P_J|\geq J$ such that
$|\delta_{|P_J|}(S\cap P_J)-\eta|<1/J$ and $\delta_{|P_J|}(A\cap P_J)>\alpha-1/J$,
which implies $\mu^j_{|P_J|}(S\cap P)=\eta$ and $\mu^j_{|P_J|}(A\cap P_J)\geq\alpha$.
On the other hand, if $\mu_{|P|}(S\cap P)=\eta$ and
$\mu^j_{|P|}(A\cap P)\geq\alpha$, then $\SD_S^j(A)\geq\alpha$ by
the definition of $\SD_S^j$ in (\ref{sds-1}).

\medskip

Part 4: If $\SD^j_S(A)=\alpha$, then $\mu^j_{|P|}(S\cap P)=\eta$ and
$\mu^j_{|P|}(A\cap P)=\alpha'\geq\alpha$ for
some $P$ with $|P|\in\n_{j'}\setminus\n_j$
by Part 3. Clearly, $\alpha=\SD^j_S(A)\geq
\SD^j_{S\cap P}(A\cap P)\geq\mu^j_{|P|}(A\cap P)=
\alpha'$ by the definition of $\SD_S^j$.
Hence $\alpha=\alpha'$.
\hfill $\blacksquare$

\medskip

The following lemma is the internal version of
an argument similar to so-called the double counting property in
the standard literature. Let $1\ll H\leq N/2$.
Roughly speaking, if $C$ is very uniformly distributed
in $[N]$ with measure $\alpha$, then for almost all $x\in [N-H]$
the measure of $C\cap (x+[H])$ inside $x+[H]$ is $\alpha$. Since
the measure $\mu_H$ is not an internal function we use $\delta_H$ instead
and require $|\delta_H(C\cap (x+[H])-\alpha|<1/J$ for some infinite $J$
instead of $\mu_H(C\cap (x+[H])=\alpha$. The lemma is stated in a
more general case with $\V_j$ being viewed as the ``standard'' universe
in $\V_{j'}$.

\begin{lemma}\label{bigD}\quad Let $N,H\in\n_{j'}\setminus\n_j$, $H\leq N/2$,
and $C\subseteq [N]$ with $\mu^j_N(C)=\SD^j(C)=\alpha\in\real_j$
for $0\leq j<j'\leq 3$. For each $n\in\n_{j'}$ let
\begin{equation}\label{bigDn}
D_{n,H,C}:=\left\{x\in [N-H]\mid
\left|\delta_H(C\cap(x+[H]))-\alpha\right|<\frac{1}{n}\right\}.
\end{equation}
Then there exists a $J\in\n_{j'}\setminus\n_j$
such that $\mu^j_{N-H}(D_{J,H,C})=1$.
\end{lemma}

Notice that $D_{n,H,C}\subseteq D_{n',H,C}$ if $n\geq n'$.

\medskip

\noindent {\it Proof}\quad Fix $N$, $H$, and $C$. The
subscripts $H$ and $C$ in $D_{n,H,C}$ will be omitted in the proof.
If $st_j(H/N)>0$, then for every $x\in [N-H]$ we have
$\mu^j_H(x+[H])=\alpha$ by the supremality of $\alpha$. Hence the maximal $J$
with $J\leq H$ such that $|\delta_H(A\cap (x+[H]))-\alpha|<1/J$
for every $x\in [N-H]$ is in $\n_{j'}\setminus\n_j$. Now $D_J=[N-H]$ works.

Assume that $st_j(H/N)=0$. So, $\mu^j_{N-H}$ and $\mu^j_N$
coincide. If $\delta_N(D_n)\approx_j 1$ for every $n\in\n_j$, then the maximal
$J$ satisfying $|\delta_N(D_J)-1|<1/J$ must be in $\n_{j'}\setminus\n_j$
by Proposition \ref{spill2}. Hence $\mu^j_N(D_J)=1$.
So we can assume that $\mu^j_N(D_n)<1$ for some $n\in\n_j$
and derive a contradiction.

Notice that for each $x\in [N-H]$, it is impossible to have
$\mu^j_H(C\cap (x+[H]))>\alpha=\SD^j(C)$ by the definition of $\SD^j$.
Let $\overline{D}_n:=[N-H]\setminus D_n$. Then $\mu^j_N(\overline{D}_n)
=1-\mu^j_N(D_n)>0$.
Notice that $x\in\overline{D}_n$ implies
$\delta_H(C\cap (x+[H]))\leq\alpha- 1/n$.
By the following double counting argument, by ignoring some
$\V_j$--infinitesimal amount inside $st_j$, we have
\begin{eqnarray*}
\lefteqn{\alpha=st_j\left(\frac{1}{H}\sum_{y=1}^H\delta_N(C-y)\right)
=st_j\left(\frac{1}{HN}\sum_{y=1}^H\sum_{x=1}^{N}\chi_C(x+y)\right)}\\
& &=st_j\left(\frac{1}{NH}\sum_{x=1}^{N}\sum_{y=1}^H
\chi_C(x+y)\right)=st_j\left(\frac{1}{N}\sum_{x=1}^{N}\delta_H(C\cap (x+[H]))\right)\\
& &=st_j\left(\frac{1}{N}\sum_{x\in D_n}\delta_H(C\cap (x+[H]))
+\frac{1}{N}\sum_{x\in\overline{D}_n}\delta_H(C\cap (x+[H]))\right)\\
& &\leq\alpha\mu^j_N(D_n)+\left(\alpha-\frac{1}{n}\right)\mu^j_N(\overline{D}_n)<\alpha
\end{eqnarray*}
which is absurd. This completes the proof.
\hfill $\blacksquare$

\medskip

Suppose $0\leq j<j'\leq 3$, $N\geq H\gg 1$ in $\n_{j'}$,
$U\subseteq [N]$, $A\subseteq S\subseteq [N]$,
$0\leq\alpha\leq\eta\leq 1$, and $x\in [N]$.  For each $n\in\n_j$ let
$\xi(x,\alpha,\eta,A,S,U,H,n)$ be the following internal statement:
\begin{equation}\label{uniform}
\begin{array}{rcl}
\vspace{0.1in}
|\delta_H(x+[H])\cap U)-1|&<&1/n,\\
\vspace{0.1in}
|\delta_H((x+[H])\cap S)-\eta|&<&1/n,\mbox{ and}\\
|\delta_H((x+[H])\cap A)-\alpha|&<&1/n.
\end{array}
\end{equation}
The statement $\xi(x,\alpha,\eta,A,S,U,H,n)$ infers
that the densities of $A,S,U$ in the interval $x+[H]$ go to
$\alpha,\eta,1$, respectively, as $n\to\infty$ in $\n_j$. The statement $\xi$
will be referred a few times in Lemma \ref{key} and its proof.

\medskip

The following lemma is the application of Lemma \ref{bigD} to the
sets $U,S,A$ simultaneously.

\medskip

\begin{lemma}\label{bigS}\quad Let $N\in\n_{j'}\setminus\n_j$,
$U\subseteq [N]$, and $A\subseteq S\subseteq [N]$ be such that
$\mu^j_N(U)=1$, $\mu^j_N(S)=\SD(S)=\eta$, and
$\mu^j_N(A)=\SD^j_S(A)=\alpha$ for some $\eta,\alpha\in\real_j$
and $0\leq j<j'\leq 3$. For any $n,h\in\n_{j'}$ let
\begin{equation}\label{bigGnh}
G_{n,h}:=\{x\in [N-h]\mid \V_{j'}\models\xi(x,\alpha,\eta,A,S,U,h,n)\}.
\end{equation}
\begin{enumerate}
\item[(a)] For each $H\in\n_{j'}\setminus\n_j$ with $H\leq N/2$
there exists a $J\in\n_{j'}\setminus\n_j$ such that $\mu^j_{N-H}(G_{J,H})=1$;
\item[(b)] For each $n\in\n_j$, there is an $h_n\in\n_j$ with $h_n>n$
such that $\delta_N(G_{n,h_n})>1-1/n$.
\end{enumerate}
\end{lemma}

\noindent {\it Proof}\quad Part (a): Applying Lemma \ref{bigD} for
$U$ and $S$ we can find $J_1,J_2\in\n_{j'}\setminus\n_j$
such that $\mu^j_{N-H}(D_{J_1,H,U})=1$ and $\mu^j_{N-H}(D_{J_2,H,S})=1$
where $D_{n,h,C}$ is defined in (\ref{bigDn}) and $\alpha$ is replaced by
 $1$ for $U$ and $\eta$ for $S$. Let $G':=D_{J_1,H,U}\cap D_{J_2,H,S}$.
For each $n\leq \min\{J_1,J_2\}$ let
\[\overline{G}''_n:=\left\{x\in [N-H]\mid
\delta_H(A\cap(x+[H]))>\alpha+\frac{1}{n}\right\},\mbox{ and}\]
\[\underline{G}''_n:=\left\{x\in [N-H]\mid
\delta_H(A\cap(x+[H]))<\alpha-\frac{1}{n}\right\}.\]
Notice that both $\overline{G}''_n$ and $\underline{G}''_n$ are $\V_{j'}$--internal.
If $\mu^j_{N-H}(\overline{G}''_n)>0$ for some $n\in\n_j$,
then $\overline{G}''_n\cap G'\not=\emptyset$.
Let $x_0\in\overline{G}''_n\cap G'$. Then we have $\mu^j_H(S\cap(x_0+[H]))=\eta$ and
$\mu^j_H(A\cap(x_0+[H]))>\alpha+1/n$, which contradicts $\SD^j_S(A)=\alpha$. Hence
$\delta_{N-H}(\overline{G}''_n)\approx_j 0$ for every $n\in\n_j$. By
Proposition \ref{spill2} we can
find $J_+\in\n_{j'}\setminus\n_j$ such that $\mu^j_{N-H}(\overline{G}''_n)=0$ for any
$n\leq J_+$. If $\mu^j_{N-H}(\underline{G}''_n)>0$ for some $n\in\n_j$,
then $\mu^j_{N-H}(\overline{G}_m'')>0$
for some $m\in\n_j$ by the fact that $\mu^j_{N-H}(A)=\alpha$.
Hence $\delta_{N-H}(\underline{G}_n'')\approx_j 0$ for every $n\in\n_j$.
By Proposition \ref{spill2} again we can find $J_-\in\n_{j'}\setminus\n_j$ such that
$\mu^j_{N-H}(\underline{G}''_n)=0$ for any $n\leq J_-$.
The proof is complete by setting $J:=\min\{J_1,J_2,J_+,J_-\}$ and
\[G_{J,H}:=(D_{J,H,U}\cap D_{J.H,S})
\setminus(\overline{G}''_{J}\cup\underline{G}''_J).\]

Part (b): Suppose Part (b) is not true. Then there exists an $n\in\n_j$ such that
$\delta_{N-h}(G_{n,h})\leq 1-1/n$
for any $h>n$ in $\n_j$. By Proposition \ref{spill2}
there is an $H\in\n_{j'}\setminus\n_j$
such that $\delta_{N-H}(G_{n,H})\leq 1-1/n$. By Part (a) there is a $J\gg n$
such that $\mu^j_{N-H}(G_{J,H})=1$. We have a contradiction because $n<J$ and hence
$G_{J,H}\subseteq G_{n,H}$.
\hfill $\blacksquare$

\medskip

Notice that for a given $n$ one can choose $h_n$ to be the least
such that $\delta_N(G_{n,h_n})>1-1/n$ in Lemma \ref{bigS} (b).
So we can assume that $h_n$ is an internal function of $n$.
Hence we can assume that $G_{n,h_n}$ is also an internal function
of $n$.

\section{Mixing Lemma}\label{RegularitySection}

\quad We work within $\V_{j'}$ for $0<j'\leq 3$ in this section.
Any unspecified sets are $\V_{j'}$--internal.
The letter $V$ will sometimes be used for a set other than
the standard/nonstandard universes in \S\ref{nsa}. Hopefully,
no confusion will arise.
The following standard lemma is a consequence of Szemer\'{e}di's Regularity Lemma
in \cite{szemeredi}. The proof of the lemma can be found in the appendix
 of \cite{tao}.

\begin{lemma}\label{regularity}
Let $U,W$ be finite sets, let
$\epsilon>0$, and for each $w\in W$, let $E_w$ be a subset of\, $U$. Then
there exists a partition $U=U_1\cup U_2\cup\cdots\cup U_{n_{\epsilon}}$ for some
$n_{\epsilon}\in\n_0$, and real numbers
$0\leq c_{u,w}\leq 1$ in $\real_0$ for $u\in[n_{\epsilon}]$ and $w\in W$
such that for any set $F\subseteq U$, one has
\[
\left||F\cap E_w|-\sum_{u=1}^{n_{\epsilon}}c_{u,w}|F\cap U_u|\right|\leq \epsilon |U|
\]
for all but $\epsilon |W|$ values of $w\in W$.
\end{lemma}

The following lemma, the nonstandard version
of so--called mixing lemma in \cite{tao},
can be derived from Lemma \ref{regularity}. We present a proof
similar to the proof in \cite{tao} in a nonstandard setting.
Part (i) and Part (ii) of the lemma are used to prove part (iii)
and only Part (iii) will be referred in the proof of Lemma \ref{key}.

\begin{lemma}[Mixing Lemma]\label{mixing}\quad Let $N\in\n_{j'}\setminus\n_0$,
$A\subseteq S\subseteq [N]$, $1\ll H\leq N/2$, and $R\subseteq [N-H]$
 be an a.p.\ with $|R|\gg 1$ such that
 \begin{equation}\label{condition1.1}
 \mu_N(S)=\SD(S)=\eta>0,\,\mu_N(A)=\SD_S(A)=\alpha>0,
 \end{equation}
 \vspace{-0.3in}
\begin{equation}\label{condition1.2}
\mu_H((x+[H])\cap S)=\eta,\,\mbox{ and }\,\mu_H((x+[H])\cap A)=\alpha
\end{equation}
for every $x\in R$. Then the following are true.
\begin{enumerate}
\item[(i)] For any set $E\subseteq [H]$ with $\mu_H(E)>0$, there is an
$x\in R$ such that
\[\mu_H(A\cap (x+E))\geq \alpha\mu_H(E);\]
\item[(ii)] Let $m\gg 1$ be such that the van der Waerden number
$\Gamma\left(3^{m},m\right)\leq |R|$.
For any internal partition $\{U_n\mid n\in [m]\}$ of $[H]$
there exists an $m$--a.p. $P\subseteq R$, a set $I\subseteq [m]$
with $\mu_H(U_I)=1$ where $U_I=\bigcup\{U_n\mid n\in I\}$,
and an infinitesimal $\epsilon>0$ such that
\[|\delta_H(A\cap (x+U_n))-\alpha\delta_H(U_n)|\leq\epsilon\delta_H(U_n)\]
 for all $n\in I$ and all $x\in P$;
\item[(iii)] Given an internal collection of sets
$\{E_w\subseteq [H]\mid w\in W\}$ with $|W|\gg 1$
and $\mu_H(E_w)>0$ for every $w\in W$,
there exists an $x\in R$ and $T\subseteq W$
such that $\mu_{|W|}(T)=1$ and
\[\mu_H(A\cap (x+E_w))=\alpha\mu_H(E_w)\]
 for every $w\in T$.
\end{enumerate}
\end{lemma}

\noindent {\it Proof}\quad Part (i):
Assume that (i) is not true. For each $x\in R$
let $r_x$ be such that $\delta_H(A\cap(E+x))=(\alpha-r_x)\delta_H(E)$.
Then $r_x$ must be positive non-infinitesimal. We can
set $r:=\min\{r_x\mid x\in R\}$ since the function $x\mapsto r_x$
is internal. Clearly, the number $r$ is positive non-infinitesimal. Hence
$\delta_H(A\cap (E+x))\leq (\alpha-r)\delta_H(E)$ for all $x\in R$.
Notice that by (\ref{condition1.1}) and (\ref{condition1.2}),
for $\mu_H$--almost all $y\in [H]$ we have
$\mu_{|R|}(S\cap(y+R))=\eta$ which implies that
for $\mu_H$--almost all $y\in [H]$
we have $\mu_{|R|}(A\cap(y+R))=\alpha$. So
\begin{eqnarray*}
\lefteqn{\alpha\mu_H(E)\approx\frac{1}{H}\sum_{y\in E}
\frac{1}{|R|}\sum_{x\in R}\chi_{A}(x+y)=\frac{1}{|R|}\sum_{x\in R}
\frac{1}{H}\sum_{y=1}^H\chi_{A\cap(E+x)}(x+y)}\\
& &\leq\frac{1}{|R|}\sum_{x\in R}(\alpha-r)\delta_H(E)
=(\alpha-r)\delta_H(E)\approx(\alpha-st(r))\mu_H(E)<\alpha\mu_H(E),
\end{eqnarray*}
which is absurd.

\medskip

Part (ii):
To make the argument explicitly internal we use $\delta_H$ instead of $\mu_H$.
For each $t\in\n_{j'}$, $x\in R$, and $n\in [m]$ let
\[c^{t}_n(x)=\left\{\begin{array}{cl}
\medskip
1 &\,\mbox{ if }\,\delta_H((x+U_{n})\cap
A)\geq\left(\alpha+\frac{1}{t}\right)\delta_H(U_n),\\
\medskip
0 &\,\mbox{ if }\,\left(\alpha-\frac{1}{t}\right)\delta_H(U_n)<
\delta_H((x+U_n)\cap
A)<\left(\alpha+\frac{1}{t}\right)\delta_H(U_n),\\
-1 &\,\mbox{ if }\,\delta_H((x+U_n)\cap
A)\leq\left(\alpha-\frac{1}{t}\right)\delta_H(U_n).
\end{array}\right.\]
and let $c^t:P\to \{-1,0,1\}^{[m]}$ be such that $c^t(x)(n)=c^t_n(x)$.
For each $t\in\n_0$, since the van der Waerden number $\Gamma(3^m,m)\leq |R|$,
there exists an $m$--a.p.\ $P_t\subseteq R$
such that $c^t(x)=c^t(x')$ for any $x,x'\in P_t$. For each $x\in P_t$ let
\[\begin{array}{rcl}
\medskip
I^+_t & = & \{n\in [m]\mid c^t(x)(n)=1\},\\
\medskip
I^-_t & = & \{n\in [m]\mid c^t(x)(n)=-1\},\mbox{ and} \\
I_t & = & [m]\setminus (I^+_t\cup I^-_t), \mbox{ and}
 \end{array}\]
 \[\begin{array}{rcl}
 \medskip
 U^+_t & = & \bigcup\{U_n\mid n\in I^+_t\},\\
 \medskip
 U^-_t & = & \bigcup\{U_n\mid n\in I^-_t\},\,\mbox{ and}\\
U_t & = & [H]\setminus (U^+_t\cup U^-_t).
\end{array}\]
Clearly, $\delta_H((x+U^-_t)\cap A)\leq(\alpha-1/t)\delta_H(U^-_t)$
because $U^-_t$ is a disjoint union of the $U_n$'s for $n\in I^-_t$.
Since $t\in\n_0$ we have that
$\mu_H(U^-_t)=0$ by (i) with $P_t$ in the
place of $R$ and $U^-_t$ in the place of $E$.
Notice that $\delta_H(A\cap(x+U^+_t))\geq
(\alpha+1/t)\delta_H(U^+_t)$. Since $\alpha\geq\mu_H(A\cap(x+U^+_t))
\geq(\alpha+1/t)\mu_H(U^+_t)$, we have that $\mu_H(U^+_t)<1$,
which implies $\mu_H(U_t)>0$.
If $\mu_H(U^+_t)>0$, then $\delta_H(A\cap(x+U^+_t))\geq
(\alpha+1/t)\delta_H(U^+_t)$ implies $\mu_H(A\cap(x+U_t))<\alpha\mu_H(U_t)$
for all $x\in P_t$, which again contradicts (i). Hence $\mu_H(U^+_t)=0$
and therefore, $\delta_H(U_t)>1-1/t$ is true for every $t\in\n_0$.

Since the set of all $t\in\n_{j'}$ with $\delta_H(U_t)>1-1/t$
is $\V_{j'}$--internal, by Proposition \ref{spill2} there is a $J\gg 1$ such that
$\delta_H(U_J)>1-1/J\approx 1$. The proof of (ii) is completed by letting
$P:=P_J$, $I:=I_J$, and $U_I:=U_J$.

\medskip

Part (iii): Choose a sufficiently large positive infinitesimal
$\epsilon$ satisfying that there is an internal partition of
$[H]=U_0\cup U_1\cup\cdots\cup U_{m}$ and
real numbers $0\leq c_{n,w}\leq 1$
for each $n\in [m]$ and $w\in W$ such that the van der Waerden number
$\Gamma(3^{m},m)\leq |R|$, and for any internal set $F\subseteq [H]$
there is a $T_F\subseteq W$ with
$|W\setminus T_F|\leq \epsilon|W|$ such that
\begin{equation}\label{approx1}
\left||F\cap E_w|-\sum_{n=1}^{m}c_{n,w}|F\cap U_n|\right|\leq \epsilon H
\end{equation}
for all $w\in T_F$. Notice that such $\epsilon$ exists because
if $\epsilon$ is a standard positive real,
then $m=n_{\epsilon}$ is in $\n_0$.
From (\ref{approx1}) with $F$ being replaced by $[H]$ we have
\begin{equation}\label{approx1-2}
\left||E_w|-\sum_{n=1}^{m}c_{n,w}|U_n|\right|\leq \epsilon H
\end{equation}
for all $w\in T_{[H]}$.
By (ii) we can find a $P\subseteq R$ of length $m$,
a positive infinitesimal $\epsilon_1$, and $I\subseteq [m]$
where, for some $x\in P$,
\[I:=\left\{n\in [m]\mid |\delta_H((x+U_n)\cap A)-\alpha\delta_H(U_n)|
<\epsilon_1\delta_H(U_n)\right\}\]
($I$ is independent of the choice of $x$),
and $V:=\bigcup\{U_n\mid n\in I\}$ with $\mu_H(V)=1$.
Let $I'=[m]\setminus I$ and $V'=[H]\setminus V$.
Then for each $w\in T:=T_{[H]}\cap T_{(A-x)\cap [H]}$ we have
\begin{eqnarray*}
\lefteqn{\left|\delta_H(A\cap (x+E_w))-\alpha\delta_H(E_w)\right|}\\
& &\leq\frac{1}{H}\left(\left||A\cap (x+E_w)|
-\sum_{n\in[m]}c_{n,w}|A\cap (x+U_n)|\right|\right.\\
& &\quad +\left|\sum_{n\in[m]}c_{n,w}|A\cap (x+U_n)|
-\sum_{n\in[m]}c_{n,w}\alpha|U_n|\right|\\
& &\quad\left. +\left|\alpha\sum_{n\in[m]}c_{n,w}|U_n|
-\alpha|E_w|\right|\right)\\
& &\leq\epsilon+\frac{1}{H}\sum_{n\in I}c_{n,w}\epsilon_1|U_n|
+2\delta_H(V')+\alpha\epsilon\\
& &\leq\epsilon+\epsilon_1\delta_H(V)+2\delta_H(V')+\alpha\epsilon\approx 0.
\end{eqnarray*}
Hence $\mu_H(A\cap (x+E_w))=\alpha\mu_H(E_w)$ for all
 $w\in T$. Notice that
$\mu_{|W|}(T)=1$ because $\epsilon\approx 0$ and
$\mu_{|W|}(T_{[H]})=\mu_{|W|}(T_{[H]\cap(A-x)})=1$.
\hfill $\blacksquare$

\medskip

The set $S$ in Lemma \ref{mixing}, although
seems unnecessary, is needed in the proof of Lemma \ref{key}.

\section{Proof of Szemer\'{e}di's Theorem}\label{allk}

\quad We work within $\V_2$ in this section except in
the proof of Claim 1 in Lemma \ref{key} where $\V_3$ is needed.

Szemer\'{e}di's theorem is an easy consequence of
Lemma \ref{key}, denoted by ${\bf L}(m)$ for all $m\in[k]$.
For an integer $n\geq 2k+1$ define an interval $C_n\subseteq [n]$ by
\begin{equation}\label{middle}
C_n:=\left[\left\lceil
\frac{kn}{2k+1}\right\rceil,\,\left\lfloor\frac{(k+1)n}{2k+1}\right\rfloor\right].
\end{equation}
The set $C_n$ is the subinterval of $[n]$
in the middle of $[n]$ with the length $\lfloor n/(2k+1)\rfloor\pm\iota$
for $\iota=0$ or $1$. If $n\gg 1$, then $\mu_n(C_n)=1/(2k+1)$.
For notational convenience we denote
\begin{equation}\label{constantD}
D:=3k^3\,\mbox{ and }\,\eta_0:=1-\frac{1}{D}.
\end{equation}

$\blacklozenge$: Fix a $K\in\n_1\setminus\n_0$. The number $K$ is the length
of an interval which will play an important role in Lemma \ref{key}.
Keeping $K$ unchanged is one of the advantages from nonstandard analysis,
which is unavailable in the standard setting.

There is a summary of ideas used in the proof of Lemma \ref{key} right after
the proof. It explains some motivation of the steps taken in the proof.

\medskip

\begin{lemma}[${\bf L}(m)$]\label{key}
Given any $\alpha>0$, $\eta>\eta_0$,
any $N\in\n_2\setminus\!\n_1$, and any
$A\subseteq S\subseteq [N]$ and $U\subseteq [N]$ with
\begin{equation}\label{condition2}
\mu_N(U)=1,\mu_N(S)=\SD(S)=\eta,
\mbox{ and }\,\mu_{N}(A)=\SD_S(A)=\alpha,
\end{equation}
the following are true:
\begin{enumerate}
\item[{\bf L}]\hspace{-0.07in}$_1(m)(\alpha,\eta,N,A,S,U,K)$:
 There exists a $k$--a.p.\
$\vec{x}\subseteq U$ with $\vec{x}\oplus [K]\subseteq [N]$ satisfying the
statement $(\forall n\in\n_0)\,\xi(\vec{x}(l),\alpha,\eta,A,S,U,K,n)$ for
$l\in[k]$, where $\xi$ is defined in (\ref{uniform}), and there exist
$T_l\subseteq C_K$ with $\mu_{|C_K|}(T_l)=1$ where $C_K$ is defined
in (\ref{middle}) and $V_l\subseteq [K]$ with $\mu_K(V_l)=1$ for every
$l\geq m$, and
collections of $k$--a.p.'s
\[\begin{array}{rcl}
\medskip
\p&:=&\bigcup\{\p_{l,t}\mid t\in T_l
\,\mbox{ and }\,l\geq m\}\,\mbox{ and}\\
\Q&:=&\bigcup\{\Q_{l,v}\mid v\in V_l\,\mbox{ and }\,l\geq m\}\,\mbox{ such that}
\end{array}\]
\begin{equation}\label{plt}
\p_{l,t}\subseteq\{p\sqsubseteq (\vec{x}\oplus[K])\cap U\mid
\forall l'<m\,(p(l')\in A)\mbox{ and } p(l)=\vec{x}(l)+t\}
\end{equation}
satisfying $\mu_K(\p_{l,t})=\alpha^{m-1}/k$ for all $l\geq m$ and $t\in T_l$, and
\begin{equation}\label{qlv}
\Q_{l,v}=\{q\sqsubseteq\vec{x}\oplus [K]\mid \forall l'<m\,(q(l')\in A)
\,\mbox{ and }\,q(l)=\vec{x}(l)+v\}
\end{equation}
satisfying $\mu_K(\Q_{l,v})\leq\alpha^{m-1}$
for all $l\geq m$ and $v\in V_l$.

\item[{\bf L}]\hspace{-0.07in}$_2(m)(\alpha,\eta,N,A,S,K$): There exist a set
$W_0\subseteq S$ of $\min\{K,\lfloor 1/D(1-\eta)\rfloor\}$--consecutive integers
where $D$ is defined in (\ref{constantD})
and a collection of $k$--a.p.'s $\R=\{r_w\mid w\in W_0\}$
 such that for each $w\in W_0$ we have
$r_w(l)\in A$ for $l< m$, $r_w(l)\in S$ for $l> m$, and $r_w(m)=w$.
\end{enumerate}
\end{lemma}

\begin{remark}
\begin{enumerate}
\item[(a)] ${\bf L}_2(m)$ is an internal statement in $\V_2$.
Both ${\bf L}_1(m)$ and ${\bf L}_2(m)$ depend on $K$. Since
$K$ is fixed throughout whole proof, it, as a parameter, may be
omitted in some expressions.

\item[(b)] If $H\gg 1$ and $T\subseteq [H]$
with $\mu_H(T)>1-\epsilon$, then\, $T$ contains $\lfloor 1/\epsilon\rfloor$
consecutive integers because otherwise we have $\mu_H(T)\leq
(\lfloor 1/\epsilon\rfloor-1)/\lfloor 1/\epsilon\rfloor$ $=
1-1/\lfloor 1/\epsilon\rfloor\leq 1-1/(1/\epsilon)=1-\epsilon$.

\item[(c)] The purpose of defining $C_K$ is that if $t\in C_K$, then the
number of $k$--a.p.'s $p\sqsubseteq\vec{x}\oplus [K]$
with $p(l)=\vec{x}(l)+t$ is guaranteed to be at least $K/(k-1)$.

\item[(d)] It is not essential to require specific constant $c=1/k$
for $\mu_K(\p_{l,t})=c\alpha^{m-1}$ in ${\bf L}_1(m)$. Just
requiring that $\mu_K(\p_{l,t})\geq c\alpha^{m-1}$ for some positive
standard real $c$ is sufficient. We use more specific expression
``$\mu_K(\p_{l,t})=\alpha^{m-1}/k$'' for notational simplicity.

\item[(e)] Some ``bad'' $k$--a.p.'s in $\p$ in ${\bf L}_1(m)$ will be
thinned out so that $\R$ in ${\bf L}_2(m)$ can be constructed
from $\p$. The collection $\Q$ is only used to prevent $\p$ from being thinned
out too much. See the proof of Lemma \ref{1implies2}.

\item[(f)] It is important to notice that in ${\bf L}_1(m)$
the collection $\p_{l,t}$ is a part of the collection at the right side of
(\ref{plt}) while the collection $\Q_{l,v}$ is equal to the collection
at the right side of (\ref{qlv}).
\end{enumerate}
\end{remark}

\begin{lemma}\label{1implies2}
${\bf L}_1(m) (\alpha,\eta,N,A,S,U)$ implies
${\bf L}_2(m) (\alpha,\eta,N,A,S)$ for any
$\alpha,\eta,N,A,S,U$ satisfying the conditions of Lemma \ref{key}.
\end{lemma}

\noindent {\it Proof}\quad Assume we have obtained the $k$--a.p.\
$\vec{x}\subseteq U$ with $\vec{x}\oplus [K]\subseteq [N]$, sets
$T_l\subseteq C_K$ and $V_l\subseteq [K]$ with
$\mu_{|C_K|}(T_l)=1$ and $\mu_K(V_l)=1$, and
collections of $k$--a.p.'s $\p$ and $\Q$
as in ${\bf L}_1(m)$.

Call a $k$--a.p.\ $p\in\p_m:=\bigcup\{\p_{m,t}\mid t\in T_m\}$
{\it good} if $p(l)\in S\cap (\vec{x}(l)+[K])$
for $l\geq m$ and {\it bad} otherwise. Let $\p_m^g\subseteq\p_m$ be
the collection of all good $k$--a.p.'s and
$\p^b_m:=\p_m\setminus\p_m^g$ be the collection of all bad $k$--a.p.'s.
Let $T_m^g:=\{p(m)-\vec{x}(m)\mid p\in\p_m^g\}$.
Then $T_m^g\subseteq T_m\cap (S-\vec{x}(m))\cap C_K$.
We show that $\mu_{|C_K|}(T_m^g)> 1-D(1-\eta)$.

Let $Q:=\{q\sqsubseteq\vec{x}\oplus [K]\mid q(l')\in A\,\mbox{ for }\,l'<m\}$.
Notice that
\[\p_m^b\subseteq\bigcup_{l\geq m}\{q\in Q\mid q(l)\not\in S\}\]
and for each $v\in V_l$,  $q\in\Q_{l,v}$
iff $q\in Q$ and $q(l)=\vec{x}(l)+v$.
\begin{eqnarray*}
\lefteqn{\hspace{-0.3in}\mbox{Hence, }\,
|\p_m^b|\leq\sum_{l=m}^k\sum_{w\in[K]\setminus(S-\vec{x}(l))}
|\{q\in Q\mid q(l)=\vec{x}(l)+w\}|}\\
& &\leq\sum_{l=m}^k\left(\sum_{w\in[K]\setminus V_l}
|\{q\in Q\mid q(l)=\vec{x}(l)+w\}|
+\sum_{v\in V_l\setminus(S-\vec{x}(l))}|\Q_{l,v}|\right)\\
& &\leq K\sum_{l=m}^k(|[K]\setminus V_l|+|V_l\setminus (S-\vec{x}(l))|\alpha^{m-1}).
\end{eqnarray*}
\vspace{-0.3in}

\[\mbox{So }\,\,|\p_m^g|=|\p_m|-|\p_m^b|
\geq \sum_{t\in T_m}|\p_{m,t}|
-K\sum_{l=m}^k(|[K]\setminus V_l|+|V_l\setminus (S-\vec{x}(l))|\alpha^{m-1}).\]
Notice that $\mu_K([K]\setminus V_l)=0$. Hence we have
$|[K]\setminus V_l|/|C_K|\approx 0$ and
\begin{eqnarray*}
\lefteqn{\mu_{|C_K|}(T_m^g)\cdot\frac{\alpha^{m-1}}{k}=st\left(\frac{1}{|C_K|}
\sum_{t\in T_m^g}\frac{1}{K}|\p_{m,t}|\right)
\geq st\left(\frac{1}{|C_K|K}|\p_m^g|\right)}\\
& &\geq st\left(\frac{1}{|C_K|}\sum_{t\in T_m}\frac{1}{K}|\p_{m,t}|
-\frac{1}{|C_K|}\sum_{l=m}^k(|V_l
\setminus (S-\vec{x}(l))|\alpha^{m-1})\right)\\
& &\geq\mu_{|C_K|}(T_m)\cdot\frac{\alpha^{m-1}}{k}
-(2k+1)k(1-\eta)\cdot\alpha^{m-1}\\
& &=\left(\frac{1}{k}-(2k+1)k(1-\eta)\right)\cdot\alpha^{m-1},
\end{eqnarray*}
which implies $\displaystyle
\mu_{|C_K|}(T_m^g)\geq 1-(2k+1)k^2(1-\eta)>1-D(1-\eta)$.
Recall that $T^g_m\subseteq C_K$. Hence $\vec{x}(m)+T^g_m$
contains a set $W_0$ of $\lfloor 1/D(1-\eta)\rfloor$ consecutive integers.
So, ${\bf L}_2(m)$ is proven if we let $\R:=\{r_w\mid w\in W_0\}$ where
$r_w$ is one of the $k$--a.p.'s in $\p_m^g$ such that $r_w(m)=\vec{x}(m)+w$.
 \hfill $\blacksquare$

\bigskip

The idea of the proof of Lemma \ref{1implies2} is due to Szemer\'{e}di. See
\cite{tao}.

\bigskip

\noindent {\it Proof of Lemma \ref{key}}\quad We prove
${\bf L}(m)$ by induction on $m$. By Lemmas \ref{1implies2}
it suffices to prove ${\bf L}_1(m)$.

For ${\bf L}(1)$, given any $\alpha>0$, $\eta>\eta_0$,
$N\in\n_2\setminus\n_1$, $A$, $S$, and $U$ satisfying
(\ref{condition2}), by Lemma \ref{bigS} (b) we can find a $k$--a.p.\
$\vec{x}\subseteq [N]$ such that
$(\forall n\in\n_0)\,\xi(\vec{x}(l),\alpha,\eta,A,S,U,K,n)$ is true
for $l\in[k]$,
where $\xi$ is defined in (\ref{uniform}).
For each $l\in [k]$ let $T_l=C_K\cap (U-\vec{x}(l))$ and $V_l=[K]$.
For each $l\in[k]$, $t\in T_l$, and $v\in V_l$ let
\[\begin{array}{rcl}
\medskip
\p_{l,t}&:=&\{p\sqsubseteq(\vec{x}\oplus[K])\cap U \mid p(l)=\vec{x}(l)+t\}\\
\Q_{l,v}&:=&\{q\sqsubseteq(\vec{x}\oplus[K]) \mid q(l)=\vec{x}(l)+v\}.
\end{array}\]
Clearly, we have $\mu_K(\p_{l,t})\geq 1/(k-1)>1/k$. By some pruning
we can assume that $\mu_K(\p_{l,t})=1/k$.
It is trivial that $\mu_K(\Q_{l,v})\leq 1$ and
$q\in\Q_{l,v}$ iff $q(l)=\vec{x}(l)+v$ for each
$q\sqsubseteq\vec{x}\oplus [K]$. This completes the proof of
${\bf L}_1(1) (\alpha,\eta,N,A,S,U)$. ${\bf L}_2(1)(\alpha,\eta,N,A,S)$ follows from
Lemma \ref{1implies2}.

\medskip

Assume ${\bf L}(m-1)$ is true for some $2\leq m\leq k$.

\medskip

We now prove ${\bf L}(m)$. Given any $\alpha>0$ and $\eta>\eta_0$,
fix $N\in\n_2\setminus\n_1$, $U\subseteq [N]$, and
$A\subseteq S\subseteq [N]$ satisfying (\ref{condition2}).
For each $n\in\n_1\setminus\n_0$, by Lemma \ref{bigS} (b),
there is an $h_n>n$ in $\n_1$
and $G_{n,h_n}\subseteq [N]$ defined in
(\ref{bigGnh}) such that $d_n:=\delta_{N-h_n}(G_{n,h_n})>1-1/n$.
Notice that $d_n\approx_1\mu^1_{N-h_n}(G_{n,h_n})>\eta_0$
because $n\gg 1$ and $\mu_{N-h_n}(G_{n,h_n}))=1$.
Let $\eta^1_n:=\mu^1_{N-h_n}(G_{n,h_n})$ and fix an $n\in\n_1\setminus\n_0$.

\medskip

{\bf Claim 1}\quad The following internal statement $\theta(n,A,N)$ is true:

$\exists W\subseteq [N]\,\exists\R\,(W\,
\mbox{ is an a.p.}\,\wedge |W|\geq \min\{K,\lfloor 1/2D(1-d_n)\rfloor\}\,\wedge\,
\R=\{r_w\mid w\in W\}$ is a collection of $k$--a.p.'s such that
\[\forall w\in W\,((\forall l\geq m)\,(r_w(l)\in\,G_{n,h_n})\,\wedge\,
r_w(m-1)=w\,\wedge\,(\forall l,l'\leq m-2)\]
\[((A\cap (r_w(l)+[h_n]))-r_w(l)
=(A\cap (r_w(l')+[h_n]))-r_w(l'))).\]

\medskip

\noindent {\it Proof of Claim 1}\quad
Working in $\V_2$ by considering $\V_1$ as the standard universe, we can find
$P\subseteq [N]$ with $|P|\in\n_2\setminus\n_1$ by Lemma \ref{equiv}
and Part 2 of Property \ref{property} such that
\[\SD^1(G_{n,h_n})=\mu^1_{|P|}(P\cap G_{n,h_n})=\SD^1(G_{n,h_n}\cap P)=\eta^1_n.\] For each
$x\in P\cap G_{n,h_n}$ let $\tau_x=((x+[h_n])\cap A)-x$. Since there are at most
$2^{h_n}\in\n_1$ different $\tau_x$'s and $|P|\gg 2^{h_n}$,
we can find one, say, $\tau_n\subseteq [h_n]$ such that the set
\[B_n:=\{x\in P\cap G_{n,h_n}\mid \tau_x=\tau_n\}\]
satisfies $\mu^1_{|P|}(B_n)\geq \eta^1_n/2^{h_n}>0$.
Notice that $\mu^0_{|P|}(B_n)$ could be $0$.

Let $P'\subseteq P$ with $|P'|=N'\in\n_2\setminus\n_1$ be
such that $\mu^1_{N'}(G_{n,h_n}\cap P')=\eta^1_n$ and
\[\begin{array}{rcl}
\vspace{0.1in}
\beta^1_n &:=&\mu^1_{N'}(B_n\cap P')=\SD^1_{G_{n,h_n}\cap P}(B_n\cap P)\\
 &=&\SD^1_{G_{n,h_n}\cap P'}(B_n\cap P')\geq
\mu^1_{|P|}(B_n)>0
\end{array}\]
by Part 4 of Lemma \ref{equiv} and Part 2 of Property \ref{property}.
Let $d$ be the common difference of the a.p.\ $P'$
and $\varphi:P'\to [N']$ be the order-preserving bijection, i.e.,
\[\varphi(x):=1+(x-\min P')/d.\] Let $B':=\varphi[B_n\cap P']$ and
$S':=\varphi[G_{n,h_n}\cap P']$.
We have that $B',S',N'$ and $\beta^1_n,\eta^1_n$ in the place of
$A,S,N$ and $\alpha,\eta$ satisfy the $\V_1$--version of (\ref{condition2})
with $\mu$, $\SD$\, and $\SD_S$ being replaced by $\mu^1$, $\SD^1$,\, and $\SD^1_{S'}$.

Let $N''=i_2(N')$, $B''=i_2(B')$, and $S''=i_2(S')$ where $i_2$ is in Part
3 of Property \ref{property}. Recall that $i_2\res\V_1$ is an identity map.
Since $N'\in\n_2\setminus\n_1$, we have
$N''\in\n_3\setminus\n_2$.
Notice also that $\mu^1_{N''}(S'')=\SD^1(S'')=\eta^1_n$
and $\mu^1_{N''}(B'')=\SD^1_{S''}(B'')=\beta^1_n$.
By the induction hypothesis that ${\bf L}(m-1)$ is true we have
\begin{eqnarray}\label{addpredicates}
\lefteqn{(\V_2;\real_0,\real_1)\models \forall\alpha,\eta\in\real_0\,\forall
N\in\n_2\setminus\n_1\,\forall A,S\subseteq [N]}\\
& &(\alpha>0\wedge\eta>\eta_0\wedge A\subseteq S\wedge\mu_N(S)=\SD(S)=\eta\wedge
\mu_N(A)=\SD_S(A)\nonumber\\
& &\to {\bf L}_2(m-1)(\alpha,\eta,N,A,S)).\nonumber
\end{eqnarray}
Since $(\V_2;\real_0,\real_1)$ and $(\V_3;\real_1,\real_2)$ are
elementarily equivalent by Part 2 of Property \ref{property} via $i_*$, we have, by
universal instantiation, that
\begin{equation}\label{addpredicates2}
(\V_3;\real_1,\real_2)\models {\bf L}_2(m-1)
(\beta^1_n,\eta^1_n,N'',B'',S'').
\end{equation}
Notice that the right side above no longer depends on $\real_1$ or $\real_2$.
So, we have
\begin{equation}\label{addpredicates3}
\V_3\models {\bf L}_2(m-1)(i_2(\beta^1_n),
i_2(\eta^1_n),i_2(N'),i_2(B'),i_2(S'))
\end{equation}
because $i_2(\beta^1_n)=\beta^1_n$ and $i_2(\eta^1_n)=\eta^1_n$.
Since $i_2$ is a bounded elementary embedding, we have
\[\V_2\models {\bf L}_2(m-1)(\beta^1_n,\eta^1_n,N',B',S'),\]
which means that there is a set $W'\subseteq [N']$ of
$\min\{K,\lfloor 1/D(1-\eta^1_n)\rfloor\}$--consecutive integers
and a collection of $k$--a.p.'s $\R'=\{r'_w\mid w\in W'\}$ such that
for every $w\in W'$ we have $r'_w(l)\in B'$ for $l<m-1$, $r'_w(m-1)=w$,
and $r'_w(l)\in S'$ for $l\geq m$.
Notice that $\varphi^{-1}[[N']]\subseteq [N]$. Let
$W=\varphi^{-1}[W']$ and $\R=\{r_w\mid w\in W\}$, where
$r_w=\varphi^{-1}[r'_{\varphi(w)}]$,
such that for each $w\in W$ we have
$r_w(l)\in\varphi^{-1}[B']\subseteq B_n$ for $l<m-1$,
$r_w(m-1)=w$, and $r_w(l)\in\varphi^{-1}[S']\subseteq G_{n,h_n}$
for $l\geq m$. If $\eta^1_n=1$, then $|W|\geq K$. If
$\eta^1_n<1$, then $2(1-d_n)>1-\eta^1_n$. Hence
$|W|\geq\min\{K,\lfloor 1/2D(1-d_n)\rfloor\}$. \hfill $\blacksquare$ (Claim 1)

\medskip

The following claim follows from Claim 1 by Proposition \ref{spill2}.

\medskip

{\bf Claim 2}\quad There exists a $J\in\n_2\setminus\n_1$ such that the
$\theta(J,A,N)$ is true, i.e.,
$\exists W\subseteq [N]\,\exists\R\,(W\,
\mbox{ is an a.p.}\,\wedge |W|\geq \min\{K,\lfloor 1/2D(1-d_J)\rfloor\}\,\wedge\,
\R=\{r_w\mid w\in W\}$ is a collection of $k$--a.p.'s such that
$\forall w\in W\,((\forall l\geq m)\,(r_w(l)\in\,G_{J,h_J})$,
$r_w(m-1)=w$, and $(\forall l,l'\leq m-2)\,((A\cap (r_w(l)+[h_J]))-r_w(l)
=(A\cap (r_w(l')+[h_J]))-r_w(l'))))$.

\medskip

For notational convenience let $W_H:=W$ and $\R_H:=\R$ be obtained in Claim 2
and rename $H:=h_J$, $S_H:=G_{J,h_J}$, $\tau_H:=(A\cap (r_w(l)+[h_J]))-r_w(l)$
for some (or any) $w\in W_H$ and $l<m-1$. Let $\{w_s\mid 1\leq s\leq |W_H|\}$
be the increasing enumeration of $W_H$.
Notice that $H\in\n_2\setminus\n_1$.
We now go back to consider $\V_0$ as our standard universe. Notice that
$\mu_{N-H}(S_H)=1$, $|W_H|\gg 1$, and
$(\forall n\in\n_0)\,\xi(x,\alpha,\eta,A,S,U,H,n)$
is true for every $x\in S_H$ where $\xi$ is defined in (\ref{uniform}).

\medskip

{\bf Claim 3}\quad
For each $s\in\n_0$ we can find an internal
$U_s\subseteq [H]$ with $\mu_H(U_s)=1$ such that for each
$y\in U_s$ and each $l\in [k]$, $r_{w_s}(l)+y\in U$ and $(\forall n\in\n_0)
\,\xi(r_{w_s}(l)+y,\alpha,\eta,A,S,U,K,n)$ is true.

\medskip

\noindent {\it Proof of Claim 3}\quad For each $l\in[k]$ we have
$\xi(r_{w_s}(l),\alpha,\eta,A,S,U,H,n)$ is true because $r_{w_s}(l)\in S_H$.
By Lemma \ref{bigS} (a), we can find a set $G_l\subseteq r_{w_s}(l)+[H]$
with $\mu_H(G_l)=1$ such that \newline
$\xi(r_{w_s}(l)+y,\alpha,\eta,A,S,U,K,n)$ is true
for every $r_{w_s}(l)+y\in G_l$. Set
\[U_s:=\bigcap_{l=1}^k((U\cap G_l)-r_{w_s}(l)).\]
Then we have $U_s\subseteq [H]$ and $\mu_H(U_s)=1$.\hfill $\blacksquare$ (Claim 3)

\medskip

Notice that $\delta_H(\bigcap_{i=1}^s U_i)>1-1/s$.
By Proposition \ref{spill2} we can find $1\ll I\leq |W_H|$ and
\[U':=\bigcap\{U_s\mid 1\leq s\leq I\}\]
such that $\delta_H(U')>1-1/I$. Hence $\mu_H(U')=1$.
Applying the induction hypothesis for
${\bf L}_1(m-1)(\alpha,1,H,\tau_H,[H],U')$, we obtain
a $k$--a.p.\ $\vec{y}\subseteq U'$ with $\vec{y}\oplus[K]\subseteq [H]$,
$T'_l\subseteq C_K\cap U'$ with $\mu_{|C_K|}(T'_l)=1$ and
$V'_l\subseteq [K]$ with $\mu_K(V'_l)=1$ for each $l\geq m-1$,
and collections of $k$--a.p.'s
\[\begin{array}{rcl}
\vspace{0.1in}
\p'&=&\bigcup\{\p'_{l,t}\mid t\in T'_l\,\mbox{ and }\,l\geq m-1\}\,\mbox{ and}\\
\Q'&=&\bigcup\{\Q'_{l,v}\mid v\in V'_l\,\mbox{ and }\,l\geq m-1\}
\end{array}\]
such that (i) for each $l\geq m-1$ and $t\in T'_l$
we have $\mu_K(\p'_{l,t})=\alpha^{m-2}/k$
and for each $p\in\p'_{l,t}$ we have $p\sqsubseteq(\vec{y}\oplus [K])\cap U'$,
$p(l')\in \tau_H$ for $l'<m-1$, $p(l)=\vec{y}(l)+t$, and (ii)
for each $l\geq m-1$ and $v\in V'_l$ we have
$\mu_K(\Q'_{l,v})\leq\alpha^{m-2}$, and for each $q\sqsubseteq\vec{y}\oplus [K]$
we have $q\in\Q'_{l,v}$ iff
 $q(l')\in \tau_H$ for every $l'< m-1$ and $q(l)=\vec{y}(l)+v$.
For each $l\geq m$, $t\in T_l$, and $v\in V_l$
let \[E_{l,t}:=\{p(m-1)\mid p\in\p'_{l,t}\}\,\mbox{ and }
\,F_{l,v}:=\{q(m-1)\mid q\in\Q'_{l,v}\}.\]
Then $E_{l,t}, F_{l,v}\subseteq\vec{y}(m-1)+[K]$, $\mu_K(E_{l,t})=\mu_K(\p'_{l,t})
=\alpha^{m-2}/k$, and $\mu_K(F_{l,v})=\mu_K(\Q'_{l,v})\leq\alpha^{m-2}$.
Since $\vec{y}\subseteq U'$ we have that for each $l\in[k]$,
$(\forall n\in\n_0)\,\xi(r_{w_s}(l)+\vec{y}(l),\alpha,\eta,A,S,U,K,n)$ is true.

Applying Part (iii) of Lemma \ref{mixing} with
$R:=\{w_s+\vec{y}(m-1)\mid 1\leq s\leq I\}$
and $H$ being replaced by $K$
we can find $s_0\in [I]$, $T_l\subseteq T'_l$
with $\mu_{|C_K|}(T_l)=1$ and $V_l\subseteq V'_l$ with $\mu_K(V_l)=1$
for each $l\geq m$ such that for each $t\in T_l$ and $v\in V_l$ we have
\begin{equation}\label{emixing}
\begin{array}{l}
\vspace{0.1in}
\mu_K((w_{s_0}+E_{l,t})\cap((w_{s_0}+\vec{y}(m-1)+[K])\cap A))\\
\qquad =\alpha\mu_K(E_{l,t})=\alpha(\alpha^{m-2}/k)
=\alpha^{m-1}/k\,\mbox{ and}
\end{array}
\end{equation}
\begin{equation}\label{fmixing}
\begin{array}{l}
\vspace{0.1in}
\mu_K((w_{s_0}+F_{l,v})\cap((w_{s_0}+\vec{y}(m-1)+[K])\cap A))\\
\qquad\qquad =\alpha\mu_K(F_{l,t})\leq\alpha\!\cdot\!\alpha^{m-2}=\alpha^{m-1}.
\end{array}
\end{equation}

Let $\vec{x}:=r_{w_{s_0}}\oplus\vec{y}$.
Clearly, we have $\vec{x}\oplus[K]\subseteq [N]$.
We also have that $\vec{x}\subseteq U$,
$\mu_K((\vec{x}(l)+[K])\cap S)=\eta$,
and $\mu_K((\vec{x}(l)+[K])\cap A)=\alpha$
because $r_{w_{s_0}}\subseteq S_H$ and
$\vec{y}\subseteq U'\subseteq U_{s_0}$. For each $l\geq m$, $t\in T_l$, and $v\in V_l$ let
\begin{eqnarray*}
\lefteqn{\p_{l,t}:=\{r_{w_{s_0}}\oplus p\mid p\in\p'_{l,t}\,\mbox{ and}}\\
& & \,p(m-1)\in
E_{l,t}\cap(((w_{s_0}+\vec{y}(m-1)+[K])\cap A)-w_{s_0})\},\\
& &\hspace{-0.3in}\Q_{l,v}:=\{r_{w_{s_0}}\oplus q\mid q\in\Q'_{l,t}\,\mbox{ and}\\
& & \,q(m-1)\in
F_{l,v}\cap(((w_{s_0}+\vec{y}(m-1)+[K])\cap A)-w_{s_0})\}.
\end{eqnarray*}
Then $\mu_K(\p_{l,t})=\alpha^{m-1}/k$ by (\ref{emixing}).
If $\bar{q}\sqsubseteq\vec{x}\oplus [K]$, then there is a $q\sqsubseteq
\vec{y}\oplus [K]$ such that $\bar{q}=r_{w_{s_0}}\oplus q$. If
$\bar{q}(l')\in A$ for $l'<m$ and $v\in V_l$
for some $l\geq m$ such that $\bar{q}(l)=\vec{x}(l)+v$, then
$q(l')\in\tau_H$ for $l'<m-1$, $v\in V'_l$, and $q(l)=\vec{y}(l)+v$,
which imply $q\in\Q'_{l,v}$ by induction hypothesis.
Hence we have $q(m-1)\in F_{l,v}$. Clearly, $\bar{q}(m-1)=w_{s_0}+q(m-1)\in A$
implies $q(m-1)\in F_{l,v}\cap (((w_{s_0}+\vec{y}(m-1)+[K])\cap A)-w_{s_0})$.
Thus we have $\bar{q}\in\Q_{l,v}$.
Clearly, $\mu_K(\Q_{l,v})\leq\alpha^{m-1}$ by (\ref{fmixing}).

Summarizing the argument above we have that for
each $r_{w_{s_0}}\oplus p\in\p_{l,t}$
\begin{itemize}
\item $r_{w_{s_0}}(l')+p(l')\in r_{w_{s_0}}(l')+\tau_H\subseteq A$ for $l'<m-1$
because $r_{w_{s_0}}(l')\in B_H$,

\item $r_{w_{s_0}}(m-1)+p(m-1)=w_{s_0}+p(m-1)$

\qquad $\in (w_{s_0}+E_{l,t})\cap
(w_{s_0}+\vec{y}(m-1)+[K])\cap A\subseteq A$,
\item $r_{w_{s_0}}(l')+p(l')\in(\vec{x}(l')+[K])\cap U\subseteq U$
for $l'\geq m$ because of $p\subseteq U'$,
\item $r_{w_{s_0}}(l)+p(l)=r_{w_{s_0}}(l)+\vec{y}(l)+t=\vec{x}(l)+t$.
\end{itemize}
For each $\bar{q}\sqsubseteq\vec{x}\oplus [K]$, $\bar{q}\in\Q_{l,v}$
iff there is a $q\sqsubseteq \vec{y}\oplus [K]$ with
$\bar{q}=r_{w_{s_0}}\oplus q$ such that
\begin{itemize}
\item $r_{w_{s_0}}(l')+q(l')\in r_{w_{s_0}}(l')+\tau_H\subseteq A$ for $l'<m-1$
because $r_{w_{s_0}}(l')\in B_H$,
\item $r_{w_{s_0}}(m-1)+q(m-1)=w_{s_0}+q(m-1)\in A$ which is equivalent to

\qquad $w_{s_0}+q(m-1)\in (w_{s_0}+F_{l,v})\cap
(w_{s_0}+\vec{y}(m-1)+[K])\cap A\subseteq A$,
\item $r_{w_{s_0}}(l)+q(l)=r_{w_{s_0}}(l)+\vec{y}(l)+v=\vec{x}(l)+v$.
\end{itemize}
This completes the proof of ${\bf L}_1(m)(\alpha,\eta,N,A,S,U)$
as well as ${\bf L}(m)$ by Lemma \ref{1implies2}.
\hfill $\blacksquare$

\bigskip

\noindent {\bf Summary of the ideas used in the proof of Lemma \ref{key}}:
\quad
We want to use ${\bf L}_2(m-1)$ to create a sequence $\R_H$ of $k$ blocks
$r_{w}\oplus[H]$ of size
$H$ that are in arithmetic progression such that (a) the set $A$
in the first $m-2$ blocks $r_{w}(l)+[H]$ for $l\leq m-2$
are identical copies of $\tau_H$ in $[H]$, (b) the initial points
$\{r_{w}(m-1)\mid w\in W_H\}$ of all $(m-1)$-st blocks form an a.p.\
 of infinite length which should
be used when applying Part (iii) of Lemma \ref{mixing},
(c) the initial points of the rest of the blocks satisfy an appropriate
version of (\ref{uniform}) for some infinite $n=J$, i.e.,
$r_w\subseteq S_H$. Then we work
inside $[H]$, using $L_1(m-1)$
to create collections of $k$--a.p.'s $\p'$ and $\Q'$ in $[H]$ instead of $[N]$.
For applying Part (iii) of Lemma \ref{mixing}
we want to make sure, if we can, that
the $(m-1)$-st terms of all $p\in\p'_{l,t}$ form a set of positive measure
and the $(m-1)$-st terms of all $q\in\Q'_{l,v}$ form
a set of positive measure. Then mixing one $r_w$ for some $w\in W_H$
with $\p'$ and $\Q'$ at
$(m-1)$-st terms yields $\p$ and $\Q$ validating ${\bf L}_1(m)$.

Unfortunately, using ${\bf L}_1(m-1)$ to create collections $\p'$ and $\Q'$ of
$k$--a.p.'s in $[H]$
cannot guarantee that the set $E_{l,t}$ of the $(m-1)$-st terms of all
$p\in\p'_{l,t}$ and the set $F_{l,v}$ of the $(m-1)$-st terms of
all $q\in\Q'_{l,v}$ have positive measures in $[H]$ (the positive measures
can be guaranteed
when $k\leq 4$ but not for $k>4$). So, instead of getting
the sets $E_{l,t}$ and $F_{l,v}$ to have positive measures
in $[H]$ we make sure that the set $E_{l,t}$
and $F_{l,v}$ have positive measures in some
subinterval $\vec{y}(m-1)+[K]$ of length $K$ in $[H]$,
where $K$ is fixed and could be much smaller than $H$.
This is achieved by requiring $\vec{y}\subseteq U'$. When we
use the mixing lemma we want to make sure that $A$ in all of these relevant
intervals $(r_{w_s}\oplus\vec{y})(m-1)+[K]$
of length $K$ has measure $\alpha$. This requirement
is again achieved by requiring $r_{w_{s}}\subseteq S_H$ and
$\vec{y}\subseteq U'$ after shrinking $W_H$ to its initial
segment $\{w_s\mid s\in [I]\}$.

Since we want $K$ to be significantly smaller than $N$
we assume that $N$ and $K$ is at least one universe apart. Hence $N$ must be
at least in $\n_2\setminus\n_1$. If we want to apply ${\bf L}_1(m-1)$
for $H$ instead of $N$, then $H$ must also be at least in
$\n_2\setminus\n_1$. If we want $B_H$ to have a $\V_2$--standard
positive measure, $N'$ must be in $\n_4$ in order to use the
$\V_4$--version of ${\bf L}(m-1)$ with $\V_2$ being considered
as the ``standard'' universe. But $N'$ can only be guaranteed
one universe apart from $H$ by Definition \ref{density} even though
$N$ is assumed to be at least two universes apart from $H$.
Therefore, we use $h_n\in\n_1$ (instead of $H$) which leads to $N'\in\n_2\setminus\n_1$
and use $i_2$ to lift $N'$ to $N''\in\n_3\setminus\n_2$ while keeping $h_n\in\n_1$.
This allows the use of $\V_3$--version of ${\bf L}_2(m-1)$ with $\V_1$ being
considered as the ``standard'' universe. Then spill $h_n$ over
to $H\in\n_2\setminus\n_1$ and apply ${\bf L}_1(m-1)$ with
$N, A, S, U$ being replaced by $H, \tau_H, [H], U'$ to obtain
desired collections
$\p'$ and $\Q'$ of $k$--a.p.'s. All these steps rely on the fact that
${\bf L}_2(m-1)$ is an internal statement with internal parameters.

\begin{theorem}[E. Szemer\'{e}di, 1975]\label{szemeredi}
Let $k\in\n_0$. If $D\subseteq\n_0$ has positive upper density,
then $D$ contains nontrivial $k$-term arithmetic progressions.
\end{theorem}

\noindent {\it Proof}\quad It suffices to find a nontrivial
$k$--a.p.\ in $i_0(D)$. Let $P$ be an a.p.\ such that $|P|\gg 1$ and
$\mu_{|P|}(i_0(D)\cap P)=\SD(D)=\alpha$. Then $\alpha>0$
because $\alpha$ is greater than or equal to the upper density of $D$.
Let $A=i_0(D)\cap P$.
Without loss of generality, we can assume $P=[N]$ for some $N\gg 1$.
We can also assume that $N\in\n_2\setminus\n_1$ because otherwise replace
$N$ by $i_1(N)$ and $A$ by $i_1(A)$.
Then we have $\mu_N(A)=\SD(A)=\alpha$. Set $U=S=[N]$. Trivially,
$\mu_N(S)=\SD(S)=\eta=1$, $A\subseteq S$, and
$\SD_S(A)=\SD(A)=\alpha$. To start with $k'=k+1$
instead of $k$, we have many nontrivial $k'$--a.p.'s
$p\in\p$ such that $p(l)\in A$ for $l\leq k'-1=k$
in ${\bf L}_1(k')$. So there must be many
nontrivial $k$--a.p.'s in $A\subseteq i_0(D)$. By $\V_0\prec\V_2$,
there must be nontrivial $k$--a.p.'s in $D$.
\hfill $\blacksquare$

\section{A Question}\label{question}

The construction of the nonstandard universes above requires the existence
of a non-principal ultrafilter $\f_0$ on $\n_0$ and a well-order $\lhd$
on some $\V(\real_0,z)$, which is a consequence
of $\mathsf{ZFC}$. Notice that $\mathsf{ZF}$ cannot guarantee
the existence of $\f_0$ although $\mathsf{ZF}$ plus the existence of $\f_0$
and $\lhd$ is strictly weaker than $\mathsf{ZFC}$.
However, assuming the existence of $\f_0$ and $\lhd$ may be avoided
by the axiomatic approach of nonstandard analysis developed in \cite{HK}.
In \cite{HK} two systems of axioms $\mathsf{SPOT}$ and $\mathsf{SCOT}$ are introduced.
Roughly speaking, $\mathsf{SPOT}$ contains $\mathsf{ZF}$ plus some primitive tools
for nonstandard analysis and $\mathsf{SCOT}$ contains $\mathsf{ZF}$ plus the axiom of
dependent choice and some primitive tools for nonstandard analysis.
It is shown in \cite{HK} that $\mathsf{SPOT}$ is a conservative extension of
$\mathsf{ZF}$ and sufficient for developing basic calculus
while $\mathsf{SCOT}$ is a conservative extension of $\mathsf{ZF}$ plus
the axiom of dependent choice and
sufficient for developing basic calculus and Lebesgue integration.

\begin{question}
Can the nonstandard proof of Szemer\'{e}di's theorem in \S\ref{allk} be carried
out in $\mathsf{SCOT}$ or even in $\mathsf{SPOT}$?
\end{question}



\section*{Acknowledgments} 
The author
would like to thank the American Institute of Mathematics
which sponsored the workshop
{\it Nonstandard Methods in Additive Combinatorics},
where he had an opportunity to attend Terence Tao's lectures
and learn from Tao's interpretation of Szemer\'{e}di's
original proof of Szemer\'{e}di's theorem \cite{tao}.
The author would also like to thank Steven Leth, Isaac Goldbring,
Mikhail Katz, Michael Benedikt, Karel Hrb\'{a}\v{c}ek, and the anonymous referee
for comments, suggestions, and correcting some mistakes and typos
in earlier versions of the paper.

\bibliographystyle{amsplain}

\begin{thebibliography}{99}
\bibitem{CK}
Chen Chung Chang and H.\ Jerome Keisler.
\newblock {\em Model Theory.}
\newblock (3rd
edition), North-Holland, 1990.

\bibitem{dinasso}
Mauro Di Nasso.
\newblock {\em Hypernatural numbers as
ultrafilters.}
\newblock in Nonstandard Analysis for the Working
Mathematician (Springer, Dordrecht, 2015), 443--474.

\bibitem{DGL}
Mauro Di Nasso, Isaac Goldbring, and Martino Lupini.
\newblock {\em Nonstandard Methods in Ramsey Theory and Combinatorial Number Theory.}
\newblock Lecture Notes in Mathematics book series, volume 2239, Springer, 2019.

\bibitem{goldblatt}
Robert Goldblatt.
\newblock {\em Lectures on the
hyperreals--an introduction to
nonstandard analysis.}
\newblock Springer, 1998.

\bibitem{HK}
Karel Hrb\'{a}\v{c}ek and Mikhail G.\ Katz.
\newblock Infinitesimal analysis without the Axiom of Choice.
\newblock {\em Annals of Pure and Applied Logic,}
\newblock 172 (6), June 2021

\bibitem{luperi}
Lorenzo Luperi Baglini.
\newblock Hyperintegers and
Nonstandard Techniques in Combinatorics of Numbers.
\newblock PhD Dissertation (2012),
\newblock University of Siena, arXiv: 1212.2049.

\bibitem{LW}
Peter A.\ Loeb and Manfred P.\ H.\ Wolff, editors.
\newblock {\em Nonstandard Analysis
for the Working Mathematician--Second edition.}
\newblock Springer, Dordrecht, 2015.

\bibitem{szemeredi}
Endre Szemer\'{e}di.
\newblock On sets of integers containing
no $k$ elements in arithmetic progression.
\newblock  Collection
of articles in memory of Juri\v{i} Vladimirovi\v{c} Linnik.
\newblock {\em Acta Arithmatica.}
\newblock 27 (1975): 199--245.

\bibitem{tao}
Terence Tao.
\newblock Szemer\'{e}di's proof of Szemer\'{e}di's theorem.
\newblock  {\em Acta Mathematica Hungarica.}
\newblock 161 (2020): 443--487.
\newblock https://terrytao.files.wordpress.com/2017/09/szemeredi-proof1.pdf

\bibitem{tao2}
Terence Tao.
\newblock A nonstandard analysis proof
of Szemer\'{e}di's theorem.\newline
\newblock https://terrytao.wordpress.com/2015/07/20/a-nonstandard-analysis-proof-of-szemeredis-theorem/

\bibitem{tao3}
Terence Tao.
\newblock {Szemeredi's proof of Szemeredi's theorem.}\newline
\newblock {\em https://terrytao.wordpress.com/2017/09/12/szemeredis-proof-of-szemeredis-theorem/}.

\end{thebibliography}


\begin{dajauthors}
\begin{authorinfo}[pgom]
  Renling Jin\\
  College of Charleston\\
  Charleston, South Carolina, USA\\

  jinr\imageat{}cofc\imagedot{}edu \\
  \url{http://jinr.people.cofc.edu}
\end{authorinfo}

\end{dajauthors}

\end{document}